\documentclass[11pt]{article}
\usepackage{amsfonts}
\usepackage{graphics}
\usepackage{indentfirst}
\usepackage{color}
\usepackage{cite}
\usepackage{latexsym}
\usepackage[paper=a4paper, left=2.1cm, right=2.1cm, top=2.2cm, bottom=1.5cm, headheight=5.5pt, footskip=0.8cm, footnotesep=0.8cm, centering, includefoot]{geometry}
\usepackage{amsmath}
\allowdisplaybreaks
\usepackage{amssymb}
\usepackage[dvips]{epsfig}
\usepackage{amscd}

\newtheorem{theorem}{Theorem}[section]
\newtheorem{remark}{Remark}[section]

\newtheorem{lemma}{Lemma}[section]
\newtheorem{corollary}{Corollary}[section]
\newtheorem{proposition}{Proposition}[section]

\DeclareMathOperator{\curl}{curl}

\DeclareMathOperator{\divv}{div}

\makeatletter
\@addtoreset{equation}{section}
\makeatother
\makeatletter
\@addtoreset{equation}{section}
\makeatother

\title{Global well-posedness to three-dimensional full compressible magnetohydrodynamic equations with vacuum
\thanks{Yang Liu is supported by National Natural Science Foundation of China (No. 11901288). Xin Zhong is supported by National Natural Science Foundation of China (No. 11901474) and the Innovation Support Program for Chongqing Overseas Returnees (No. cx2019130).
}
}

\author{Yang Liu\thanks{College of Mathematics, Changchun Normal
University, Changchun 130032, People's Republic of China ({\tt liuyang0405@ccsfu.edu.cn}). }
 \quad Xin Zhong\thanks{Corresponding author. School of Mathematics and Statistics, Southwest University, Chongqing 400715,  People's Republic of China ({\tt xzhong1014@amss.ac.cn}).
 }
 }

\date{ }

\begin{document}
\maketitle

\begin{abstract}
This paper studies the Cauchy problem for three-dimensional viscous,
compressible, and heat conducting magnetohydrodynamic equations with vacuum
as far field density. We prove the global existence and uniqueness of strong solutions provided that the quantity $\|\rho_0\|_{L^\infty}+\|b_0\|_{L^3}$ is suitably small and the viscosity coefficients satisfy $3\mu>\lambda$. Here, the initial velocity and initial temperature could be large. The
assumption on the initial density do not exclude that the initial density may vanish in a subset of $\mathbb{R}^3$ and that it can be of a nontrivially compact support. Our result is an extension of the works of Fan and Yu \cite{FY09} and Li et al. \cite{LXZ13}, where the local strong solutions in three dimensions and the global strong solutions for isentropic case were obtained, respectively. The analysis is based on some new mathematical techniques and some new useful energy estimates. This paper can be viewed as the first result concerning the global existence of strong solutions with vacuum at infinity in some classes of large data in higher dimension.
\end{abstract}


Math Subject Classification: 35Q55; 76N10; 76W05

\section{Introduction}
Let $\Omega\subset\mathbb{R}^3$ be a domain, the motion of a viscous, compressible, and heat conducting magnetohydrodynamic (MHD) flow in $\Omega$ can be described by full compressible MHD equations (see \cite[Chapter 3]{LQ2012}):
\begin{equation}
\left\{
\begin{array}{ll}
\displaystyle
 \rho_t+\divv(\rho u)=0,\\[3pt]
\rho u_t+\rho u\cdot\nabla u-\mu\Delta u-(\lambda+\mu)\nabla\divv u+\nabla p=\curl b\times b,\\[3pt]
c_v\rho(\theta_t+u\cdot\nabla\theta)+p\divv u-\kappa\Delta\theta=\mathcal{Q}(\nabla u)
+\nu|\curl b|^2,\\[3pt]
b_t-b\cdot\nabla u+u\cdot\nabla b+b\divv u=\nu\Delta b,\\[3pt]
\divv b=0,
\end{array}
\right.\label{a1}
\end{equation}
where the unknowns $\rho\ge 0$, $u\in\Bbb R^3$, $\theta\ge 0$, and $b\in\Bbb R^3$ are
the density, velocity, pressure, absolute temperature, and magnetic field,
respectively; $p=R\rho\theta$, with positive constant $R$, is the pressure, and
\begin{align}
\mathcal{Q}(\nabla u)=\frac{\mu}{2}|\nabla u+(\nabla u)^\top|^2
+\lambda(\divv u)^2,
\end{align}
with $(\nabla u)^\top$ being the transpose of $\nabla u$. The constant viscosity coefficients $\mu$ and $\lambda$ satisfy the physical restrictions
\begin{align}\label{1.2}
\mu>0, \quad 2\mu+3\lambda\ge 0.
\end{align}
Positive constants $c_\nu$, $\kappa$, and $\nu$ are the heat capacity, the ratio of the heat conductivity coefficient over the heat capacity, and the magnetic diffusive coefficient, respectively.

Let $\Omega=\mathbb{R}^3$ and we consider the Cauchy problem of \eqref{a1} with $(\rho,u,\theta,b)$ vanishing at
infinity (in some weak sense) with given initial data $\rho_0$, $u_0$, $\theta_0$, and $b_0$, as
\begin{align}\label{a2}
(\rho, u, \theta, b)|_{t=0}=(\rho_0, u_0, \theta_0, b_0), \quad x\in\Bbb R^3.
\end{align}

The compressible MHD equations govern the motion of electrically conducting fluids such as plasmas, liquid metals, and electrolytes. They consist of a coupled system of compressible Navier-Stokes equations of fluid dynamics and Maxwell's equations of electromagnetism. Besides their wide physical applicability (see e.g., \cite{D2017}), the MHD system are also of great interest in mathematics. As a coupled system, the issues of well-posedness and dynamical behaviors of compressible MHD equations are rather complicated to investigate because of the strong coupling and interplay interaction between the fluid motion and the magnetic field. Their distinctive features make analytic studies a great challenge but offer new opportunities.
Furthermore, the differences in behaviors of solutions between isentropic and non-isentropic fluid flows are believed to be significant (see \cite{DF2006,LQ2012,HW08,HW10}).

On the one hand, for isentropic case, Suen and Hoff \cite{SH12} proved the global-in-time existence of weak solutions in three space dimensions with initial data small in $L^2$ and initial density positive and essentially bounded. As emphasized in many related papers (refer to \cite{HS1991,HS2001,LXY1998,X1998,XY2013} for instance), the possible appearance of vacuum produces new difficulty in mathematical analysis, so it is interesting to study the solutions with vacuum. Hu and Wang \cite{HW10} showed the global weak solutions with vacuum with large initial data in terms of the Lions' compactness framework of renormalized solutions \cite{L1998}. The global-in-time weak solutions for a non-resistive fluid in two dimensions were obtained recently in \cite{LS2019}.
Moreover, for the global well-posedness of strong solutions with vacuum, Li et al. \cite{LXZ13} and L{\"u} et al. \cite{LSX16} established the global existence and uniqueness of strong solutions to the 3D case and 2D case, respectively, provided the smooth initial data are of small total energy, which generalize similar results for strong solutions of the isentropic compressible Navier-Stokes equations obtained by Huang et al. \cite{HLX2012} and Li and Xin \cite{LX2019}, respectively. Later, by removing the crucial assumption that the initial total energy is small, Hong et al. \cite{HHPZ} improved the result of \cite{LXZ13} and proved the global classical strongs as long as the adiabatic exponent is close to 1 and $\nu$ is suitably large.

On the other hand, for non-isentropic case \eqref{a1}, Kawashima \cite{K83} first obtained the global existence and uniqueness of classical solutions in multi-dimension when the initial data are close to a non-vacuum equilibrium in $H^3$-norm (see also \cite{PG13}). Using the entropy method, Ducomet and Feireisl \cite{DF2006} studied the the global existence of weak solutions by introducing the entropy equation rather than the thermal equation \eqref{a1}$_3$. Meanwhile, Hu and Wang \cite{HW08} considered global-in-time weak solutions of \eqref{a1} instead of the entropy equation used in \cite{DF2006}. Non-uniqueness of global-in-time weak solutions for an inviscid fluid in two dimensions was investigated in \cite{FL2020}. For local well-posedness of strong solutions with vacuum, Fan and Yu \cite{FY09} established the local existence and uniqueness of strong solutions to \eqref{a1}--\eqref{a2}. Zhong \cite{Z20202} investigated the 2D case of \eqref{a1} with $\kappa=\nu=0$ via weighted energy method.
However, to the best of our knowledge, global well-posedness theory for strong solutions with vacuum to \eqref{a1} in multi-dimension cannot be available. In fact, the main aim of this paper is to deal with the global existence and uniqueness of strong solutions to the 3D Cauchy problem \eqref{a1}--\eqref{a2} in some homogeneous Sobolev spaces with vacuum at infinity for the density and the temperature.

Before formulating our main result, we first explain the notations and conventions used throughout this paper.
 For simplicity, in what follows, we denote
\begin{align*}
\int_{\Bbb R^3}fdx=\int fdx, \quad c_v=\kappa=R=\nu=1.
\end{align*}
For $1\le p\le \infty$ and integer $k\ge 0$, the standard homogeneous and inhomogeneous
Sobolev spaces as follows:
\begin{align*}
\left\{
\begin{array}{ll}
\displaystyle
L^p=L^p(\Bbb R^3),~ ~ W^{k, p}=L^p\cap D^{k, p},~~ H^k=W^{k, 2}£¬\\[3pt]
  D^{k, p}=\{u\in L_{loc}^1(\Bbb R^3): \|\nabla^ku\|_{L^p}<\infty\}, ~D^k=D^{k, 2},\\[3pt]
  D_0^1=\{u\in L^6(\Bbb R^3): \|\nabla u\|_{L^2}<\infty\}.
\end{array}
\right.
\end{align*}

Let $E_0$ be the specific energy defined by
\begin{align*}
E_0=\frac{|u_0|^2}{2}+\theta_0.
\end{align*}

Our main result can be stated as follows.
\begin{theorem}\label{thm1}
Let $3\mu>\lambda$. For given numbers $K>0$ (which may be arbitrarily large), $q\in (3, 6)$, and $\bar{\rho}>0$, assume that the initial data $(\rho_0\ge 0, u_0, \theta_0\ge 0, b_0)$ satisfies
\begin{equation}
\left\{
\begin{array}{ll}
\displaystyle
\rho_0\le \bar{\rho}, ~\rho_0\in L^1\cap H^1\cap W^{1, q}, ~(u_0, \theta_0)\in D_0^1\cap D^{2, 2}, \\[3pt]
 \sqrt{\rho_0}E_0+\sqrt{\rho_0}u_0\in L^2,\ b_0\in H^2,\
\divv b_0=0,\\[3pt]
 \|\sqrt{\rho_0}u_0\|_{L^2}^2+\|\nabla u_0\|_{L^2}^2+\|\sqrt{\rho_0}E_0\|_{L^2}^2+\|b_0\|_{H^1}^2=K,\\[3pt]
 \bar{\rho}+\|b_0\|_{L^3}=M_0^2,
\end{array}
\right.\label{qq}
\end{equation}
and the compatibility conditions
\begin{align}\label{1.6}
\begin{cases}
-\mu\Delta u_0-(\mu+\lambda)\nabla\divv u_0+\nabla(\rho_0\theta_0)
-\curl b\times b=\sqrt{\rho_0}g_1,\\[3pt]
\Delta\theta_0+\mathcal{Q}(\nabla u_0)+|\curl  b|^2=\sqrt{\rho_0}g_2,
\end{cases}
\end{align}
with $g_1, g_2\in L^2$.
There exists a small positive
constant $\epsilon_0$ depending only on  $\mu$, $\lambda$, $\|\rho_0\|_{L^1}$, and $K$
 such that if
\begin{align*}
M_0\le \epsilon_0,
\end{align*}
then the problem \eqref{a1}--\eqref{a2} has a unique global strong solution $(\rho\ge 0, u, \theta\ge 0, b)$ satisfying
\begin{equation}
\left\{
\begin{array}{ll}
\displaystyle
\rho \in C([0, \infty); L^1\cap H^1\cap W^{1, q}), \ \rho_t\in C([0, \infty); L^2\cap L^q),\\
(u, b, \theta)\in C([0, \infty); D_0^1\cap D^2)\cap L_{loc}^2([0, \infty); D^{2, q}),\ b\in C([0, \infty); H^2),\\
(u_t, b_t, \theta_t)\in L_{loc}^2([0, \infty); D_0^1), ~(b_t, \sqrt{\rho}u_t, \sqrt{\rho}\theta_t)\in L_{loc}^\infty([0, \infty); L^2).
\end{array}
\right.\label{trr}
\end{equation}
\end{theorem}

\begin{remark}
Theorem \ref{thm1} is the first result concerning the global existence of
strong solutions to the full compressible magnetohydrodynamic equations with vacuum in spatial multi-dimension. Moreover, the conclusion in Theorem \ref{thm1} generalizes the theory of isentropic case in Li et al \cite{LXZ13} to the non-isentropic case. In particular, the initial energy is allowed
to be large when $\|\rho_0\|_{L^\infty}+\|b_0\|_{L^3}$ is suitably small.
\end{remark}

\begin{remark}
It should be noted that there is no need to require any smallness condition on the initial velocity $u_0$ and initial temperature $\theta_0$
 for the global existence of solutions.
\end{remark}

\begin{remark}
It is very interesting to investigate the global existence and uniqueness of strong solutions to the initial boundary value problem of \eqref{a1} under various boundary conditions for $(u,\theta,b)$. Some new ideas are needed to handle these cases. This will be left for future studies.
\end{remark}

If $b\equiv b_0\equiv0$, Theorem \ref{thm1} directly yields the following global existence theorem for the full compressible Navier-Stokes equations.
\begin{theorem}\label{thm2}
Let $3\mu>\lambda$. For given numbers $K>0$ (which may be arbitrarily large), $q\in (3, 6)$, and $\bar{\rho}>0$, assume that the initial data $(\rho_0, u_0, \theta_0\ge 0)$ satisfies
\begin{equation}
\left\{
\begin{array}{ll}
\displaystyle
0\le \rho_0\le \bar{\rho}, ~\rho_0\in L^1\cap H^1\cap W^{1, q}, \\[3pt]
 \sqrt{\rho_0}E_0+\sqrt{\rho_0}u_0\in L^2,\
 (u_0, \theta_0)\in D_0^1\cap D^{2, 2}, \\[3pt]
 \|\sqrt{\rho_0}u_0\|_{L^2}^2+\|\nabla u_0\|_{L^2}^2+\|\sqrt{\rho_0}E_0\|_{L^2}^2=K,
\end{array}
\right.\label{zqq}
\end{equation}
and the compatibility conditions
\begin{align}
\begin{cases}
-\mu\Delta u_0-(\mu+\lambda)\nabla\divv u_0+\nabla(\rho_0\theta_0)
=\sqrt{\rho_0}g_1,\\[3pt]
\Delta\theta_0+\mathcal{Q}(\nabla u_0)=\sqrt{\rho_0}g_2,
\end{cases}
\end{align}
with $g_1, g_2\in L^2$.
There exists a small positive
constant $\epsilon_0$ depending only on  $\mu$, $\lambda$, $\|\rho_0\|_{L^1}$, and $K$
 such that if
\begin{align*}
\bar{\rho}\le \epsilon_0,
\end{align*}
then the problem \eqref{a1}--\eqref{a2}  with $b\equiv0$ has a unique global strong solution $(\rho\ge 0, u, \theta\ge 0)$ satisfying
\begin{equation}
\left\{
\begin{array}{ll}
\displaystyle
\rho \in C([0, \infty); L^1\cap H^1\cap W^{1, q}), \ \rho_t\in C([0, \infty); L^2\cap L^q),\\
(u, \theta)\in C([0, \infty); D_0^1\cap D^2)\cap L_{loc}^2([0, \infty); D^{2, q}),\\
(u_t, \theta_t)\in L_{loc}^2([0, \infty); D_0^1),
~(\sqrt{\rho}u_t, \sqrt{\rho}\theta_t)\in L_{loc}^\infty([0, \infty); L^2).
\end{array}
\right.\label{ztrr}
\end{equation}
\end{theorem}

\begin{remark}
Since the assumption $3\mu>\lambda$ is weaker than $2\mu>\lambda$ due to $\mu>0$, Theorem \ref{thm2} extends the result in Li \cite{JL19} where the  global existence of strong solution was established provided that $\|\rho_0\|_{L^\infty}(\|\rho_0\|_{L^3}
+\|\rho_0\|_{L^\infty}^2\|\sqrt{\rho_0}u_0\|_{L^2}^2)
(\|\nabla u_0\|_{L^2}^2+\|\rho_0\|_{L^\infty}\|\sqrt{\rho_0}E_0\|_{L^2}^2)$ is sufficiently small and $2\mu>\lambda$.
\end{remark}

\begin{remark}
We note that in \cite{HL18}, Huang and Li studied the Cauchy problem of full
compressible Navier-Stokes equations in $\mathbb{R}^3$, and they obtained the existence and uniqueness of global classical solutions provided that the
initial energy is small. However, the initial density and initial temperature are not allowed to vanish at infinity. Such assumptions play a crucial role for some estimates in \cite{HL18}.
\end{remark}

We now make some comments on the analysis for Theorem \ref{thm1}. To prove the global existence of strong solutions, we establish a crucial proposition (Proposition \ref{p1}) which implies that the terms in Serrin-type criterion (see Lemma \ref{pyy}) will never blow up in finite time when $3\mu>\lambda$ and the initial data is small in some sense (refer to Section 4 for more details). This together with the contradiction arguments indicates that the strong solution exists globally in time. This is the main ingredient of the proof. Compared to the isentropic case \cite{LXZ13}, due to $(\rho(x, t),b(x,t),\theta(x, t))\rightarrow (0,0,0)$ as $|x|\rightarrow\infty$, the basic energy inequality only provides us
\begin{align*}
\int\big(\rho|u|^2+|b|^2+2\rho\theta\big)dx
=\int\big(\rho_0|u_0|^2+|b_0|^2+2\rho_0\theta_0\big)dx,
\end{align*}
and there is no any useful dissipation estimate on $u$ and $b$. To overcome this difficulty, inspired by \cite{JL19,WZ17}, where the authors obtained dissipative estimate on $u$ for the full Navier-Stokes equations by using $L^3$-norm of the density and the conservation of mass, respectively, we recover the crucial dissipation estimate of the form $\int_0^T(\mu\|\nabla u\|_{L^2}+\|\nabla b\|_{L^2})dt$ in terms of $L^\infty(0,T;L^\infty)$-norm of $\rho$ (see Lemma \ref{lem32}).
Moreover, as stated in many papers (see \cite{SH12,DF2006,HW08} for example), compared with compressible Navier-Stokes equations, the presence of magnetic
field effects results in some new difficulties. To this end, we try to deal with the strong coupling term $u\cdot \nabla b$ and the strong nonlinear term $\curl b\times b$ by introducing the spatial $L^\infty(0,T;L^3)$-norm of $b$. These motivate us to impose the smallness condition
on $\|\rho_0\|_{L^\infty}+\|b_0\|_{L^3}$ to get the bound of $\|\rho\|_{L^\infty}+\|b\|_{L^3}$.
Furthermore, we tackle higher order estimates with the help of the effective viscous flux $F=(2\mu+\lambda){\rm div}u-p-\frac{1}{2}|b|^2$ (see Lemma \ref{l35}) and the upper bound of the density is obtained via commutator estimate (see Lemma \ref{lem38}). Finally, it enables
us to get $L^\infty(0,T;L^3)$ estimate of $b$ from the induction equation \eqref{a1}$_4$ and Kato-type inequality (see Lemma \ref{l39}). Combining these estimates altogether yields the desired energy-like estimate, provided that the initial data is suitably small (see Corollary \ref{qww}).

The rest of the paper is organized as follows. In Section 2, we recall some known facts and elementary inequalities which will be used later. Section 3 is devoted to the global a priori estimates. The proof of Theorem \ref{thm1} will be done in Section 4.

\section{Preliminaries}
In this section, we collect some known results and elementary inequalities which will be used later.

First, the following local existence and uniqueness of strong solutions has been established in \cite{FY09}.
\begin{lemma}\label{l22}
Assume that $(\rho_0, u_0, \theta_0, b_0)$ satisfies \eqref{qq} and \eqref{1.6}. Then there exists a small time $T>0$ and a
unique strong solution $(\rho, u, \theta, b)$ to the problem \eqref{a1}--\eqref{a2} on $\Bbb R^3\times(0, T)$.
\end{lemma}

Next, the following well-known Gagliardo-Nirenberg inequality (see \cite[Theorem]{HL}) will be used later frequently.
\begin{lemma}\label{l21}
Let $u$ belong to $L^q(\mathbb{R}^n)$ and its derivatives of order $m, \nabla^m u$, belong to $L^r(\mathbb{R}^n)$, $1\leq q, r\leq \infty$. Then
for the drivatives $\nabla^j u, 0\leq j< m$, the following inequality holds.
\begin{align}\label{bb1}
\|\nabla^j u\|_{L^p(\mathbb{R}^n)}\leq C\|\nabla^m u\|_{L^r(\mathbb{R}^n)}^\alpha\|u\|_{L^q(\mathbb{R}^n)}^{1-\alpha},
\end{align}
where
\begin{align}
\frac{1}{p}=\frac{j}{n}+\alpha\left(\frac{1}{r}-\frac{m}{n}\right)
+(1-\alpha)\frac{1}{q}
\end{align}
for all $\alpha$ in the interval
\begin{align}
\frac{j}{m}\leq \alpha\leq 1
\end{align}
(the constant $C$ depends only on $n, m, j, q, r, \alpha$), with the following exceptional cases:

(1) If $j=0, rm<n$ and $q=\infty$, then we take the additional assumption that either $u$
tends to zero at infinity or $u\in L^{\tilde{q}}(\mathbb{R}^n)$ for some finite $\tilde{q}>0$.

(2) If $1<r<\infty$, and $m-j-\frac{n}{r}$ is a nonegative integer, then
\eqref{bb1} holds only for $\alpha$ satisfying $\frac{j}{m}\leq \alpha<1$.
\end{lemma}

Finally, the following Serrin-type blow-up criterion (see \cite{HL13}) will be used to prove the global existence of the strong solution to \eqref{a1}--\eqref{a2} (see Section 4 for details).
\begin{lemma}\label{pyy}
Let the initial data  $(\rho_0, u_0, \theta_0, b_0)$ satisfy conditions in Theorem \ref{thm1}. If $T^{*}<\infty$ is the maximal time of existence for that solution obtained in Lemma \ref{l22}, then we have
\begin{align*}
\lim_{T\rightarrow T^*}(\|\rho\|_{L^\infty(0, T; L^\infty)}+\|u\|_{L^{s}(0, T; L^{r})})=\infty,
\end{align*}
with $r$ and $s$ satisfying
\begin{align*}
\frac{2}{s}+\frac{3}{r}\le 1, \ s>1, \ 3<r\le \infty.
\end{align*}
\end{lemma}

\section{A Priori Estimates}
This section is devoted to deriving the following a priori estimates for the solutions to the
Cauchy problem \eqref{a1}--\eqref{a2}. For simplicity, we denote
\begin{align*}
\psi_T:=\sup_{0\le t\le T}\big(\|\sqrt{\rho}u\|_{L^2}^2+\|\nabla u\|_{L^2}^2+\|\sqrt{\rho}E\|_{L^2}^2+\|b\|_{H^1}^2\big).
\end{align*}

\begin{proposition}\label{p1}
Assume $3\mu>\lambda$, and let the conditions in Theorem \ref{thm1} be in force. There exists a positive constant $\epsilon_0$ depending only on
$\mu$, $\lambda$, $\|\rho_0\|_{L^1}$, and $K$, such that if
\begin{align}\label{ee1}
\sup_{0\le t\le T}\|\rho\|_{L^\infty}\le 2\bar{\rho}, \ \psi_T\le 2\hbar K, \ \sup_{0\le t\le T}\|b\|_{L^3}\le 2M_0,
\end{align}
then one has
\begin{align}\label{ee2}
\sup_{0\le t\le T}\|\rho\|_{L^\infty}\le \frac{3}{2}\bar{\rho}, \ \psi_T\le \frac{7}{4}\hbar K, \ \sup_{0\le t\le T}\|b\|_{L^3}\le \frac{3}{2}M_0,
\end{align}
provided that $M_0\le \epsilon_0$. Here, the constant $\hbar=\frac{16\mu+9\lambda}{\mu}$.

\end{proposition}

The proof of Proposition \ref{p1} will be done by a series of lemmas below. For simplicity, we will use the
conventions that $C$ and $C_i\ (i = 1, 2,\cdots)$ denote various positive constants, which may depend on $\mu$, $\lambda$, $\|\rho_0\|_{L^1}$, and $K$,
 but are independent of $T$ and $M_0$.

We begin with the following lemma concerning the mass is conserved for all time, which could be found in \cite[Lemma 3.1]{WZ17}, and so we omit the detail of proof.

\begin{lemma}\label{l31}
Under the conditions of Proposition \ref{p1}, it holds that
\begin{align}\label{abb}
\int\rho dx=\int\rho_0dx.
\end{align}
\end{lemma}

\begin{lemma}\label{lem32}
Under the conditions of Proposition \ref{p1}, it holds that
\begin{align}\label{bb}
&\sup_{0\le t\le T}(\|\sqrt{\rho}u\|_{L^2}^2+\|b\|_{L^2}^2)+\int_0^T\big(\mu\|\nabla u\|_{L^2}^2+\|\nabla b\|_{L^2}^2\big)dt\nonumber\\
&\le \|\sqrt{\rho_0}u_0\|_{L^2}^2+\|b_0\|_{L^2}^2+CM_0^\frac{8}{3}\int_0^T\|\nabla\theta\|_{L^2}^2dt.
\end{align}
\end{lemma}
{\it Proof.}
Multiplying $\eqref{a1}_2$ by $u$, $\eqref{a1}_4$ by $b$, respectively,
then adding the two resulting equations together, and integrating over $\mathbb{R}^3$, and noting that
$\mu+\lambda>0$\footnote{From \eqref{1.2} and $3\mu>\lambda$, we have $5\mu+2\lambda>0$. Then by \eqref{1.2} again one gets $7\mu+5\lambda>0$, which combined with \eqref{1.2} again implies $9\mu+8\lambda>0$. This together with \eqref{1.2} once more gives $11\mu+11\lambda>0$. Thus the result follows.}, we obtain from \eqref{abb} that
\begin{align*}
& \frac12\frac{d}{dt}(\|\sqrt{\rho}u\|_{L^2}^2+\|b\|_{L^2}^2)+\mu\|\nabla u\|_{L^2}^2+(\mu+\lambda)\|\divv u\|_{L^2}^2+\|\nabla b\|_{L^2}\nonumber\\
& =\int p\divv udx
\le \|\rho\|_{L^3}\|\theta\|_{L^6}\|\divv u\|_{L^2} \notag \\
& \le (\mu+\lambda)\|\divv u\|_{L^2}^2
+C\|\rho\|_{L^3}^2\|\nabla\theta\|_{L^2}^2\nonumber\\
&\le (\mu+\lambda)\|\divv u\|_{L^2}^2
+C\|\rho\|_{L^\infty}^\frac{4}{3}\|\rho\|_{L^1}^\frac{2}{3}\|\nabla\theta\|_{L^2}^2\nonumber\\
&\le (\mu+\lambda)\|\divv u\|_{L^2}^2
+CM_0^\frac{8}{3}\|\nabla\theta\|_{L^2}^2,
\end{align*}
which implies that
\begin{align}\label{p11}
\frac{d}{dt}(\|\sqrt{\rho}u\|_{L^2}^2+\|b\|_{L^2}^2)+\mu\|\nabla u\|_{L^2}^2+\|\nabla b\|_{L^2}^2
\le CM_0^\frac{8}{3}\|\nabla\theta\|_{L^2}^2.
\end{align}
Hence the desired \eqref{bb} follows from \eqref{p11} integrated in $t$.
\hfill $\Box$

\begin{lemma}
Under the conditions of Proposition \ref{p1}, it holds that
\begin{align}\label{vvq}
\sup_{0\le t\le T}\|\sqrt{\rho}E\|_{L^2}^2+\int_0^T\|\nabla\theta\|_{L^2}^2dt
& \le \|\sqrt{\rho_0}E_0\|_{L^2}^2+\frac53\int_0^T\||u||\nabla u|\|_{L^2}^2dt+CM_0^\frac{16}{3}\int_0^T\|\nabla\theta\|_{L^2}^2dt
 \nonumber\\
&\quad +CM_0^2\int_0^T\|\nabla^2 b\|_{L^2}^2dt+CM_0^4\int_0^T\|\nabla u\|_{L^2}^2dt,
\end{align}
where $E=\frac{|u|^2}{2}+\theta$.
\end{lemma}
{\it Proof.}
For $E=\frac{|u|^2}{2}+\theta$, we infer from \eqref{a1} that
\begin{align}\label{b5}
\rho (E_t+ u\cdot\nabla E)+\divv(up)-\Delta\theta=\divv(\mathcal{S}\cdot u)+\curl b\times b
+|\curl b|^2,
\end{align}
where $\mathcal{S}=\mu(\nabla u+(\nabla u)^\top)+\lambda\divv u\Bbb I_3$ with $\mathbb{I}_3$ being the identity matrix of order $3$. Multiplying \eqref{b5} by $E$ and integrating the resultant
over $\Bbb R^3$, it follows from integration by parts and Young's inequality that
\begin{align}\label{b7}
\frac{1}{2}\frac{d}{dt}\|\sqrt{\rho}E\|_{L^2}^2+\|\nabla\theta\|_{L^2}^2
&\le -\frac{1}{2}\int\nabla\theta\cdot\nabla|u|^2dx+\int(up-\mathcal{S}\cdot u)\cdot \nabla Edx\nonumber\\
&\quad+C\int(|u||b|^2|\nabla E|+|\nabla u||b|^2E)dx
+\int |\curl b|^2Edx\nonumber\\
&\le \frac{1}{6}\|\nabla\theta\|_{L^2}^2+\frac{3}{8}\||u||\nabla u|\|_{L^2}^2+C\int\rho^2\theta^2|u|^2dx\nonumber\\
&\quad+C\int(|u||b|^2|\nabla E|+|\nabla u||b|^2E)dx+C\int|\nabla E||\nabla b||b|dx\nonumber\\
&\quad+C\int|E||\nabla^2 b||b|dx=:\sum_{i=1}^6I_i.
\end{align}
Using H\"older's, the Sobolev, and the Cauchy inequalities, we have
\begin{align}
I_3&\le C\|\sqrt{\rho}\theta\|_{L^2}\|\theta\|_{L^6}\||u|^2\|_{L^6}\|\rho\|_{L^9}^\frac{3}{2}\nonumber\\
&\le C\|\sqrt{\rho}\theta\|_{L^2}\|\nabla\theta\|_{L^2}\||u||\nabla u|\|_{L^2}\|\rho\|_{L^1}^\frac{1}{6}\|\rho\|_{L^\infty}^\frac{4}{3}\nonumber\\
&\le \frac{1}{8}\||u||\nabla u|\|_{L^2}^2+C\bar{\rho}^\frac{8}{3}\|\rho\|_{L^1}^\frac{1}{3}\|\sqrt{\rho}\theta\|_{L^2}^2\|\nabla\theta\|_{L^2}^2\nonumber\\
&\le \frac{1}{8}\||u||\nabla u|\|_{L^2}^2+CM_0^\frac{16}{3}\|\nabla\theta\|_{L^2}^2,\label{xss}\\
I_4&\le C\|u\|_{L^6}\||b|^2\|_{L^3}\|\nabla E\|_{L^2}+C\|\nabla u\|_{L^2}\||b|^2\|_{L^3}\|E\|_{L^6}\nonumber\\
&\le C\|\nabla u\|_{L^2}\|b\|_{L^\infty}\|b\|_{L^3}\|\nabla E\|_{L^2}\nonumber\\
&\le C\|b\|_{L^3}^\frac{4}{3}\|\nabla^2 b\|_{L^2}^\frac{2}{3}\|\nabla u\|_{L^2}\|\nabla E\|_{L^2}\nonumber\\
&\le \frac{1}{6}\|\nabla E\|_{L^2}^2
+C\|b\|_{L^3}^\frac{8}{3}\|\nabla u\|_{L^2}^2\|\nabla^2 b\|_{L^2}^\frac{4}{3}\nonumber\\
&\le \frac{1}{6}\|\nabla\theta\|_{L^2}^2+\frac{1}{6}\||u||\nabla u|\|_{L^2}^2
+CM_0^2\|\nabla^2 b\|_{L^2}^2+CM_0^4\|\nabla u\|_{L^2}^6\nonumber\\
&\le \frac{1}{6}\|\nabla\theta\|_{L^2}^2+\frac{1}{6}\||u||\nabla u|\|_{L^2}^2
+CM_0^2\|\nabla^2 b\|_{L^2}^2+CM_0^4\|\nabla u\|_{L^2}^2,\\
I_5+I_6&\le C\|E\|_{L^6}\|\nabla^2 b\|_{L^2}\|b\|_{L^3}+
C\|\nabla E\|_{L^2}\|\nabla b\|_{L^6}\|b\|_{L^3}\nonumber\\
&\le C\|\nabla E\|_{L^2}\|b\|_{L^3}\|\nabla^2 b\|_{L^2}\nonumber\\
&\le \frac{1}{6}\|\nabla\theta\|_{L^2}^2+\frac{1}{6}\||u||\nabla u|\|_{L^2}^2+CM_0^2\|\nabla^2 b\|_{L^2}^2. \label{zxss}
\end{align}
Inserting \eqref{xss}--\eqref{zxss} into \eqref{b7} gives rise to
\begin{align}\label{b8}
\frac{d}{dt}\|\sqrt{\rho}E\|_{L^2}^2+\|\nabla\theta\|_{L^2}^2
\le \frac53\||u||\nabla u|\|_{L^2}^2
+CM_0^\frac{16}{3}\|\nabla\theta\|_{L^2}^2
+CM_0^2\|\nabla^2 d\|_{L^2}^2+CM_0^4\|\nabla u\|_{L^2}^2.
\end{align}
Then integrating \eqref{b8} in $t$ leads to the desired \eqref{vvq}.
\hfill $\Box$

Next, motivated by \cite{Z19}, we can improve the basic estimate obtained in Lemma \ref{lem32}.
\begin{lemma}\label{lem34}
Under the conditions of Proposition \ref{p1}, there exists a positive constant $c_1$ depending on $\mu$, $\lambda$, $\|\rho_0\|_{L^1}$, and $K$,
 but independent of $T$ and $M_0$, such that
\begin{align}\label{sk}
\sup_{0\le t\le T}\|\rho^\frac{1}{4}u\|_{L^4}^4+c_1\int_0^T\||u||\nabla u|\|_{L^2}^2dt
&\le CM_0^\frac{5}{3}\|\nabla u_0\|_{L^2}^4
+CM_0^\frac{16}{3}\int_0^T
\|\nabla\theta\|_{L^2}^2dt
\nonumber\\
&\quad
+CM_0^2\int_0^T\|\nabla^2 b\|_{L^2}^2dt+CM_0^4\int_0^T\|\nabla u\|_{L^2}^2dt.
\end{align}
\end{lemma}
{\it Proof.}
Multiplying $\eqref{a1}_2$ by $4|u|^2u$ and integrating the resulting equation over $\Bbb R^3$ yield
\begin{align}\label{bbc}
&\frac{d}{dt}\int\rho|u|^4dx+4\int|u|^2\big(\mu|\nabla u|^2+(\mu+\lambda)|{\rm div}u|^2+2\mu\big|\nabla|u|\big|^2\big)dx\nonumber\\
&\le 4\int{\rm div}(|u|^2u)pdx-8(\lambda+\mu)\int{\rm div}u|u|u\cdot\nabla|u|dx+C\int|u|^2|\nabla u||b|^2dx.
\end{align}
For the last term of the right-hand side of \eqref{bbc}, one obtains from H\"older's and Gagliardo-Nirenberg
inequalities that, for any $\eta_1\in (0, 1)$,
\begin{align*}
C\int|u|^2|\nabla u||b|^2dx&\le 4\mu\eta_1\int|u|^2|\nabla u|^2dx+C(\eta_1)\int|u|^2|b|^4dx\nonumber\\
&\le 4\mu\eta_1\int|u|^2|\nabla u|^2dx+C(\eta_1)\|u\|_{L^6}^2\|b\|_{L^\infty}^2\|b\|_{L^3}^2\nonumber\\
&\le 4\mu\eta_1\int|u|^2|\nabla u|^2dx+C(\eta_1)\|\nabla u\|_{L^2}^2\|b\|_{L^3}^\frac{8}{3}\|\nabla^2 b\|_{L^2}^\frac{4}{3}\nonumber\\
&\le 4\mu\eta_1\int|u|^2|\nabla u|^2dx+C\|b\|_{L^3}^2\|\nabla^2 b\|_{L^2}^2+C\|b\|_{L^3}^4\|\nabla u\|_{L^2}^6,
\end{align*}
which together with \eqref{bbc} leads to
\begin{align}
&\frac{d}{dt}\int\rho|u|^4dx+4\int|u|^2\big(\mu(1-\eta_1)|\nabla u|^2+(\mu+\lambda)|{\rm div}u|^2+2\mu\big|\nabla|u|\big|^2\big)dx\nonumber\\
&\le 4\int{\rm div}(|u|^2u)pdx-8(\lambda+\mu)\int{\rm div}u|u|u\cdot\nabla|u|dx
+C\|b\|_{L^3}^2\|\nabla^2 b\|_{L^2}^2
+C\|b\|_{L^3}^4\|\nabla u\|_{L^2}^6.
\end{align}
Consequently, we arrive at
\begin{align}\label{ssq}
&\frac{d}{dt}\int\rho|u|^4dx+4\int_{\Bbb R^3\cap\{|u|>0\}}\Big[\mu(1-\eta_1)|u|^2|\nabla u|^2+(\mu+\lambda)|u|^2|{\rm div}u|^2+2\mu|u|^2\big|\nabla|u|\big|^2\Big]dx\nonumber\\
&\le 4\int_{\Bbb R^3\cap\{|u|>0\}}{\rm div}(|u|^2u)pdx-8(\lambda+\mu)\int_{\Bbb R^3\cap\{|u|>0\}}{\rm div}u|u|u\cdot\nabla|u|dx
+C\|b\|_{L^3}^2\|\nabla^2 b\|_{L^2}^2\nonumber\\[3pt]
&\quad+C\|b\|_{L^3}^4\|\nabla u\|_{L^2}^6.
\end{align}
Direct calculations give that for $x\in \Bbb R^3\cap\{|u|>0\}$,
\begin{align}
&|u|^2|\nabla u|^2=|u|^4\Big|\nabla\Big(\frac{u}{|u|}\Big)\Big|^2+|u|^2\big|\nabla|u|\big|^2,\label{yu}\\
&|u|{\rm div}u=|u|^2{\rm div}(\frac{u}{|u|})+u\cdot\nabla|u|.
\end{align}
For $\eta_1, \eta_2\in (0, 1)$, we now define a nonnegative function as follows:
\begin{align}\label{3.19}
\phi(\eta_1, \eta_2)=\left\{
\begin{array}{ll}
\displaystyle
\frac{\mu\eta_2(3-\eta_1)}{\lambda+\eta_1\mu}, &{\rm if}~\lambda+\eta_1\mu>0,\\
0, &{\rm otherwise}.
\end{array}
\right.
\end{align}

We prove \eqref{sk} in two cases.

{\it Case 1:} we assume that
\begin{align}\label{3.20}
\int_{\Bbb R^3\cap\{|u|>0\}}|u|^4\Big|\nabla\Big(\frac{u}{|u|}\Big)\Big|^2dx\le \phi(\eta_1, \eta_2)
\int_{\Bbb R^3\cap\{|u|>0\}}|u|^2\big|\nabla|u|\big|^2dx.
\end{align}
It follows from \eqref{ssq} that
\begin{align}\label{s00}
&\frac{d}{dt}\int\rho|u|^4dx+4\int_{\Bbb R^3\cap\{|u|>0\}}Gdx\nonumber\\
&\le 4\int_{\Bbb R^3\cap\{|u|>0\}}{\rm div}(|u|^2u)pdx
+C\|b\|_{L^3}^2\|\nabla^2 b\|_{L^2}^2+C\|b\|_{L^3}^4\|\nabla u\|_{L^2}^6,
\end{align}
where
\begin{align*}
G=\mu(1-\eta_1)|u|^2|\nabla u|^2+(\mu+\lambda)|u|^2|{\rm div}u|^2+2\mu|u|^2\big|\nabla|u|\big|^2
+2(\lambda+\mu){\rm div}u|u|u\cdot\nabla|u|.
\end{align*}
To let $\int_{\Bbb R^3\cap\{|u|>0\}}Gdx$  become a good term, we shall consider $G$ first. It follows from \eqref{yu} that
\begin{align}\label{3.22}
G&=\mu(1-\eta_1)|u|^2|\nabla u|^2+(\mu+\lambda)|u|^2|{\rm div}u|^2+2\mu|u|^2|\nabla|u||^2\nonumber\\
&\quad+2(\lambda+\mu)|u|^2{\rm div}\Big(\frac{u}{|u|}\Big)u\cdot\nabla|u|+2(\lambda+\mu)|u\cdot\nabla|u||^2\nonumber\\
&=\mu(1-\eta_1)\Big(|u|^4\Big|\nabla\Big(\frac{u}{|u|}\Big)\Big|^2+|u|^2\big|\nabla|u|\big|^2\Big)
+(\lambda+\mu)\Big(|u|^2{\rm div}\Big(\frac{u}{|u|}\Big)+u\cdot\nabla|u|\Big)^2\nonumber\\
&\quad+2\mu|u|^2|\nabla|u||^2+2(\lambda+\mu)|u|^2{\rm div}\Big(\frac{u}{|u|}\Big)u\cdot\nabla|u|
+2(\lambda+\mu)|u\cdot\nabla|u||^2\nonumber\\
&= \mu(1-\eta_1)|u|^4\Big|\nabla\Big(\frac{u}{|u|}\Big)\Big|^2
+\mu(3-\eta_1)|u|^2|\nabla|u||^2-\frac{\lambda+\mu}{3}|u|^4\Big|{\rm div}\Big(\frac{u}{|u|}\Big)\Big|^2\nonumber\\
&\quad+3(\lambda+\mu)\Big(\frac{2}{3}|u|^2{\rm div}\Big(\frac{u}{|u|}\Big)+u\cdot\nabla|u|\Big)^2\nonumber\\
&\ge-(\lambda+\eta_1\mu)|u|^4\Big|\nabla\Big(\frac{u}{|u|}\Big)\Big|^2+\mu(3-\eta_1)|u|^2|\nabla|u||^2.
\end{align}
Here we have used the following facts
\begin{align*}
&(3\mu-\lambda)+4(2\mu+3\lambda)=11(\mu+\lambda)>0;\\
&\Big|{\rm div}\Big(\frac{u}{|u|}\Big)\Big|^2\le 3\Big|\nabla\Big(\frac{u}{|u|}\Big)\Big|^2.
\end{align*}
Thus, we obtain from \eqref{3.22} and \eqref{3.20} that
\begin{align}\label{3.23}
\int_{\Bbb R^3\cap\{|u|>0\}}Gdx&\ge \Big[-(\lambda+\eta_1\mu)\phi(\eta_1, \eta_2)+4\mu(3-\eta_1)\Big]
\int_{\Bbb R^3\cap\{|u|>0\}}|u|^2|\nabla|u||^2dx\nonumber\\
&\ge \mu(3-\eta_1)(1-\eta_2)\int_{\Bbb R^3\cap\{|u|>0\}}|u|^2|\nabla|u||^2dx.
\end{align}
Inserting \eqref{3.23} into \eqref{s00}, we have
\begin{align}\label{3.24}
&\frac{d}{dt}\int\rho|u|^4dx+4\mu(3-\eta_1)(1-\eta_2)\int_{\Bbb R^3\cap\{|u|>0\}}|u|^2|\nabla|u||^2dx\nonumber\\
&\le 4\int_{\Bbb R^3\cap\{|u|>0\}}{\rm div}(|u|^2u)pdx
+C\|b\|_{L^3}^2\|\nabla^2b\|_{L^2}^2+C\|b\|_{L^3}^4\|\nabla u\|_{L^2}^6\nonumber\\
&\le C\int_{\Bbb R^3\cap\{|u|>0\}}|u|^2|\nabla u|\rho\theta dx
+C\|b\|_{L^3}^2\|\nabla^2 b\|_{L^2}^2+C\|b\|_{L^3}^4\|\nabla u\|_{L^2}^6\nonumber\\
&\le \eta\int_{\Bbb R^3\cap\{|u|>0\}}|u|^2|\nabla u|^2dx
+C\int_{\Bbb R^3\cap\{|u|>0\}}\rho^2\theta^2|u|^2dx
+C\|b\|_{L^3}^2\|\nabla^3 d\|_{L^2}^2+C\|b\|_{L^3}^4\|\nabla u\|_{L^2}^6\nonumber\\
&\le \eta'\||u||\nabla u|\|_{L^2}^2+C\bar{\rho}^\frac{8}{3}\|\rho\|_{L^1}^\frac{1}{3}
\|\sqrt{\rho}\theta\|_{L^2}^2\|\nabla\theta\|_{L^2}^2
+C\|b\|_{L^3}^2\|\nabla^2 b\|_{L^2}^2+C\|b\|_{L^3}^4\|\nabla u\|_{L^2}^6\nonumber\\
&\le \eta'(1+\phi(\eta_1,\eta_2))
\int_{\Bbb R^3\cap\{|u|>0\}}|u|^2\big|\nabla|u|\big|^2dx+C\bar{\rho}^\frac{8}{3}\|\rho\|_{L^1}^\frac{1}{3}
\|\sqrt{\rho}\theta\|_{L^2}^2\|\nabla\theta\|_{L^2}^2\nonumber\\
&\quad+C\|b\|_{L^3}^2\|\nabla^2 b\|_{L^2}^2+C\|b\|_{L^3}^4\|\nabla u\|_{L^2}^6.
\end{align}
Taking $\eta'=\frac{2\mu(3-\eta_1)(1-\eta_2)}{1+\phi(\eta_1, \eta_2)}$, then we infer from \eqref{3.24} that
\begin{align}\label{fg}
&\frac{d}{dt}\int\rho|u|^4dx+2\mu(3-\eta_1)(1-\eta_2)\int_{\Bbb R^3\cap\{|u|>0\}}|u|^2|\nabla|u||^2dx\nonumber\\
&\le C\bar{\rho}^\frac{8}{3}\|\rho\|_{L^1}^\frac{1}{3}
\|\sqrt{\rho}\theta\|_{L^2}^2\|\nabla\theta\|_{L^2}^2
+C\|b\|_{L^3}^2\|\nabla^2 b\|_{L^2}^2+C\|b\|_{L^3}^4\|\nabla u\|_{L^2}^6.
\end{align}

{\it Case 2:} we assume that
\begin{align}\label{ss}
\int_{\Bbb R^3\cap\{|u|>0\}}|u|^4\Big|\nabla\Big(\frac{u}{|u|}\Big)\Big|^2dx> \phi(\eta_1, \eta_2)
\int_{\Bbb R^3\cap\{|u|>0\}}|u|^2\big|\nabla|u|\big|^2dx.
\end{align}
It follows from \eqref{bbc} that
\begin{align}
&\frac{d}{dt}\int\rho|u|^4dx+4\int\big(\mu|u|^2|\nabla u|^2+(\mu+\lambda)|u|^2|{\rm div}u|^2+2\mu|u|^2\big|\nabla|u|\big|^2\big)dx\nonumber\\
&\le 4\int{\rm div}(|u|^2u)pdx-8(\lambda+\mu)\int{\rm div}u|u|u\cdot\nabla|u|dx
+2\mu\eta_1\int_{\Bbb R^3\cap\{u>0\}}|u|^2|\nabla u|^2dx\nonumber\\
&\quad+C\|b\|_{L^3}^2\|\nabla^2 b\|_{L^2}^2
+C\|b\|_{L^3}^4\|\nabla u\|_{L^2}^6\nonumber\\
&\le C\int_{\Bbb R^3\cap\{u>0\}}p|u|^2|\nabla u|dx+4(\mu+\lambda)\int_{\Bbb R^3\cap\{u>0\}}|u|^2|\nabla|u||^2dx
+2\mu\eta_1\int_{\Bbb R^3\cap\{u>0\}}|u|^2|\nabla u|^2dx\nonumber\\
&\quad+4(\mu+\lambda)\int_{\Bbb R^3\cap\{u>0\}}|u|^2|{\rm div}u|^2dx+C\|b\|_{L^3}^2\|\nabla^2 b\|_{L^2}^2
+C\|b\|_{L^3}^4\|\nabla u\|_{L^2}^6\nonumber\\
&\le C\int_{\Bbb R^3\cap\{u>0\}}p|u|^2|\nabla|u||dx+C\int_{\Bbb R^3\cap\{u>0\}}p|u|^3\Big|\nabla\Big(\frac{u}{|u|}\Big)\Big|
+2\mu\eta_1\int_{\Bbb R^3\cap\{u>0\}}|u|^2|\nabla u|^2dx\nonumber\\
&\quad+4(\mu+\lambda)\int_{\Bbb R^3\cap\{u>0\}}|u|^2|\nabla|u||^2dx+4(\mu+\lambda)\int_{\Bbb R^3\cap\{u>0\}}|u|^2|{\rm div}u|^2dx\nonumber\\
&\quad+C\|b\|_{L^3}^2\|\nabla^2 b\|_{L^2}^2
+C\|b\|_{L^3}^4\|\nabla u\|_{L^2}^6\nonumber\\
&\le C\int_{\Bbb R^3\cap\{u>0\}}p|u|^2|\nabla|u||dx+4\mu(1-\eta_1)\eta_3\int_{\Bbb R^3\cap\{u>0\}}|u|^4\Big|\nabla\Big(\frac{u}{|u|}\Big)\Big|^2dx\nonumber\\
&\quad+4(\mu+\lambda)\int_{\Bbb R^3\cap\{u>0\}}|u|^2|\nabla|u||^2dx+4(\mu+\lambda)\int_{\Bbb R^3\cap\{u>0\}}|u|^2|{\rm div}u|^2dx\nonumber\\
&\quad+C(\eta_1, \eta_3)\int_{\Bbb R^3\cap\{u>0\}}\rho^2\theta^2|u|^2dx+C\|b\|_{L^3}^2\|\nabla^2 b\|_{L^2}^2
+C\|b\|_{L^3}^4\|\nabla u\|_{L^2}^6\nonumber\\
&\quad++2\mu\eta_1\int_{\Bbb R^3\cap\{u>0\}}|u|^2|\nabla u|^2dx\nonumber\\
&\le 4\mu\eta_1\int_{\Bbb R^3\cap\{u>0\}}|u|^2|\nabla u|^2dx+4\mu(1-\eta_1)\eta_3\int_{\Bbb R^3\cap\{u>0\}}|u|^4\Big|\nabla\Big(\frac{u}{|u|}\Big)\Big|^2dx\nonumber\\
&\quad+4(\mu+\lambda)\int_{\Bbb R^3\cap\{u>0\}}|u|^2|\nabla|u||^2dx+4(\mu+\lambda)\int_{\Bbb R^3\cap\{u>0\}}|u|^2|{\rm div}u|^2dx\nonumber\\
&\quad+C\int_{\Bbb R^3\cap\{u>0\}}p|u|^2|\nabla|u||dx+C\bar{\rho}^\frac{8}{3}\|\rho\|_{L^1}^\frac{1}{3}
\|\sqrt{\rho}\theta\|_{L^2}^2\|\nabla\theta\|_{L^2}^2\nonumber\\
&\quad+C\|b\|_{L^3}^2\|\nabla^2 b\|_{L^2}^2+C\|b\|_{L^3}^4\|\nabla u\|_{L^2}^6,
\end{align}
which together with \eqref{yu} and \eqref{ss} yields
\begin{align}\label{kk0}
&\frac{d}{dt}\int\rho|u|^4dx+f(\eta_1, \eta_2, \eta_3, \eta_4)\int_{\Bbb R^3\cap\{u>0\}}|u|^2\big|\nabla|u|\big|^2dx\nonumber\\
&\quad+4\mu(1-\eta_1)(1-\eta_3)\eta_4\int_{\Bbb R^3\cap\{u>0\}}|u|^4\Big|\nabla\Big(\frac{u}{|u|}\Big)\Big|^2dx\nonumber\\
&\le C\int_{\Bbb R^3\cap\{u>0\}}p|u|^2|\nabla|u||dx+C\bar{\rho}^\frac{8}{3}\|\rho\|_{L^1}^\frac{1}{3}
\|\sqrt{\rho}\theta\|_{L^2}^2\|\nabla\theta\|_{L^2}^2\nonumber\\
&\quad+C\|b\|_{L^3}^2\|\nabla^2 b\|_{L^2}^2+C\|b\|_{L^3}^4\|\nabla u\|_{L^2}^6.
\end{align}
where
\begin{align}
f(\eta_1, \eta_2, \eta_3, \eta_4)=4\mu(1-\eta_1)(1-\eta_3)(1-\eta_4)\phi(\eta_1, \eta_2)+8\mu-4(\lambda+\eta_1\mu),
\end{align}
 for $\eta_i\in (0, 1)\ (i=1, 2, 3, 4)$ to be decided later.

({\it Sub-case $\mathrm{1}_1$})
If $\lambda<0$, take $\eta_1=-\frac{\lambda}{m\mu}\in (0, 1)$, with the positive integer $m$ large enough, then we have
\begin{align}
\eta_1\mu+\lambda=\frac{m-1}{m}\lambda<0,
\end{align}
which combined with \eqref{3.19} implies $\phi(\eta_1, \eta_2)=0$, and hence
\begin{align}
f(\eta_1, \eta_2, \eta_3, \eta_4)=8\mu-4(\lambda+\eta_1\mu)>8\mu>0.
\end{align}

({\it Sub-case $\mathrm{1}_2$})
If $\lambda=0$, then $\phi(\eta_1, \eta_2)=\frac{\eta_2(3-\eta_1)}{\eta_1}$, and thus
\begin{align}
f(\eta_1, \eta_2, \eta_3, \eta_4)=\frac{4\mu(1-\eta_1)(1-\eta_3)(1-\eta_4)(3-\eta_1)\eta_2}{\eta_1}+8\mu-4\eta_1\mu>4\mu>0.
\end{align}

({\it Sub-case $\mathrm{1}_3$})
If $3\mu>\lambda>0$, then we have
\begin{align}
f(\eta_1, \eta_2, \eta_3, \eta_4)=\frac{4\mu^2(1-\eta_1)(1-\eta_3)(1-\eta_4)(3-\eta_1)\eta_2}{\lambda+\eta_1\mu}
+8\mu-4(\lambda+\eta_1\mu).
\end{align}
Since $f(\eta_1, \eta_2, \eta_3, \eta_4)$ is continuous w.r.t. $(\eta_1, \eta_2, \eta_3, \eta_4)$ over
$[0, 1]\times[0, 1]\times[0, 1]\times[0, 1]$, and
\begin{align}
f(0, 1, 0, 0)=\frac{12\mu^2}{\lambda}+8\mu-4\lambda>0,
\end{align}
there exists some
$(\eta_1, \eta_2, \eta_3, \eta_4)\in (0, 1)\times(0, 1)\times(0, 1)\times(0, 1)$ such that
\begin{align}
f(\eta_1, \eta_2, \eta_3, \eta_4)>0.
\end{align}
By \eqref{kk0}, Cauchy-Schwarz inequality, and H\"older's inequality, we have
\begin{align*}
&\frac{d}{dt}\int\rho|u|^4dx+f(\eta_1, \eta_2, \eta_3, \eta_4)\int_{\Bbb R^3\cap\{u>0\}}|u|^2\big|\nabla|u|\big|^2dx\nonumber\\
&\quad+4\mu(1-\eta_1)(1-\eta_3)\eta_4\int_{\Bbb R^3\cap\{u>0\}}|u|^4\Big|\nabla\Big(\frac{u}{|u|}\Big)\Big|^2dx\nonumber\\
&\le \frac{f(\eta_1, \eta_2, \eta_3, \eta_4)}{2}\int_{\Bbb R^3\cap\{u>0\}}|u|^2\big|\nabla|u|\big|^2dx+C\bar{\rho}^\frac{8}{3}\|\rho\|_{L^1}^\frac{1}{3}
\|\sqrt{\rho}\theta\|_{L^2}^2\|\nabla\theta\|_{L^2}^2\nonumber\\
&\quad+C\|b\|_{L^3}^2\|\nabla^2 b\|_{L^2}^2+C\|b\|_{L^3}^4\|\nabla u\|_{L^2}^6,
\end{align*}
that is,
\begin{align}\label{fr}
&\frac{d}{dt}\int\rho|u|^4dx+f(\eta_1, \eta_2, \eta_3, \eta_4)\int_{\Bbb R^3\cap\{u>0\}}|u|^2\big|\nabla|u|\big|^2dx\nonumber\\
&\quad+4\mu(1-\eta_1)(1-\eta_3)\eta_4\int_{\Bbb R^3\cap\{u>0\}}|u|^4\Big|\nabla\Big(\frac{u}{|u|}\Big)\Big|^2dx\nonumber\\
&\le C\bar{\rho}^\frac{8}{3}\|\rho\|_{L^1}^\frac{1}{3}
\|\sqrt{\rho}\theta\|_{L^2}^2\|\nabla\theta\|_{L^2}^2
+C\|b\|_{L^3}^2\|\nabla^2 b\|_{L^2}^2+C\|b\|_{L^3}^4\|\nabla u\|_{L^2}^6.
\end{align}
From \eqref{fg}, \eqref{fr}, and \eqref{yu}, for {\it Case 1} and {\it Case 2}, we conclude
that if $3\mu>\lambda$, there exists a constant $c_1$ such that
\begin{align}
&\sup_{0\le t\le T}\|\rho^\frac{1}{4}u\|_{L^4}^4+c_1\int_0^T\||u||\nabla u|\|_{L^2}^2dt\nonumber\\
&\le \|\rho_0^\frac{1}{4} u_0\|_{L^4}^4+C\bar{\rho}^\frac{8}{3}\int_0^T
\|\sqrt{\rho}\theta\|_{L^2}^2\|\nabla\theta\|_{L^2}^2dt
+C\int_0^T\|b\|_{L^3}^2\|\nabla^2 b\|_{L^2}^2dt
+C\int_0^T|b\|_{L^3}^4\|\nabla u\|_{L^2}^6dt\nonumber\\
&\le CM_0^\frac{2}{3}\|\nabla u_0\|_{L^2}^4
+CM_0^\frac{16}{3}\int_0^T
\|\nabla\theta\|_{L^2}^2dt+CM_0^2\int_0^T\|\nabla^2 b\|_{L^2}^2dt
+CM_0^4\int_0^T\|\nabla u\|_{L^2}^2dt.
\end{align}
Here we have used the following fact
\begin{align*}
\int\rho_0|u_0|^4dx&\le \|\rho_0\|_{L^\infty}^\frac{1}{2}\|\sqrt{\rho_0}u_0\|_{L^2}\|u_0\|_{L^6}^3
\le C\|\rho_0\|_{L^\infty}^\frac{5}{6}\|\rho_0\|_{L^1}^\frac{2}{3}\|\nabla u_0\|_{L^2}^4\le CM_0^\frac{5}{3}\|\nabla u_0\|_{L^2}^4.
\end{align*}
The proof of Lemma \ref{lem34} is completed.
\hfill $\Box$

\begin{lemma}\label{l35}
Under the conditions of Proposition \ref{p1}, it holds that
\begin{align}\label{vxx}
&\sup_{0\le t\le T}\big(\|\nabla u\|_{L^2}^2+\|\nabla b\|_{L^2}^2\big)+\int_0^T\Big(\frac{2}{\mu}\|\sqrt{\rho}u_t\|_{L^2}^2
+\|b_t\|_{L^2}^2+\|\nabla^2 b\|_{L^2}^2\Big)dt\nonumber\\
&\le \frac{15\mu+9\lambda}{\mu}\|\nabla u_0\|_{L^2}^2+2\|\nabla b_0\|_{L^2}^2
+\frac{6M_0^2}{\mu(2\mu+\lambda)}\|\sqrt{\rho_0}\theta_0\|_{L^2}^2
+CM_0^2\|\nabla b_0\|_{L^2}^2\nonumber\\
&\quad+CM_0\|\sqrt{\rho}\theta\|_{L^2}^2
+CM_0^\frac{1}{4}\int_0^T\||u||\nabla u|\|_{L^2}^2dt
+CM_0^3\int_0^T\|\nabla u\|_{L^2}^2dt,
\end{align}
provided $M_0\le \epsilon_2=\min\Big\{\epsilon_1, \Big(\frac{1}{4C_2}\Big)^\frac{5}{7},
\Big(\frac{1}{4C_3}\Big)^2\Big\}$.
\end{lemma}
{\it Proof.}
Multiplying $\eqref{a1}_2$ by $u_t$ and integrating resultant over $\Bbb R^3$, we get from integration by parts that
\begin{align}\label{mmy}
&\frac{1}{2}\frac{d}{dt}\big(\mu\|\nabla u\|_{L^2}^2+(\mu+\lambda)\|{\rm div}u\|_{L^2}^2\big)
+\|\sqrt{\rho}u_t\|_{L^2}^2\nonumber\\
&=\frac{d}{dt}\int \Big(\frac{1}{2}|b|^2{\rm div}u-b\cdot\nabla u\cdot b+p{\rm div}u\Big)dx
-\int p_t{\rm div}udx\nonumber\\
&\quad+\int(b_t\cdot\nabla u\cdot b+b\cdot\nabla u\cdot b_t-b\cdot b_t{\rm div}u)dx
-\int\rho u\cdot\nabla u\cdot u_tdx\nonumber\\
&=\frac{d}{dt}\int \Big(\frac{1}{2}|b|^2{\rm div}u-b\cdot\nabla u\cdot b+p{\rm div}u\Big)dx
-\frac{1}{2(2\mu+\lambda)}\frac{d}{dt}\int p^2dx\nonumber\\
&\quad+\int(b_t\cdot\nabla u\cdot b+b\cdot\nabla u\cdot b_t-b\cdot b_t{\rm div}u)dx
-\frac{1}{2\mu+\lambda}\int p_tFdx\nonumber\\
&\quad-\frac{1}{2(2\mu+\lambda)}\int p_t|b|^2dx-\int\rho u\cdot\nabla u\cdot u_tdx=:\sum_{i=1}^6J_i,
\end{align}
where $F=(2\mu+\lambda){\rm div}u-p-\frac{1}{2}|b|^2$.

By \eqref{ee1} and Gagliardo-Nirenberg inequality, we have
\begin{align*}
J_3&\le C\|b\|_{L^\infty}\|b_t\|_{L^2}\|\nabla u\|_{L^2}\nonumber\\
&\le C\|b\|_{L^3}^\frac{1}{3}\|\nabla^2 b\|_{L^2}^\frac{2}{3}\|b_t\|_{L^2}\|\nabla u\|_{L^2}\nonumber\\
&\le \frac{1}{2}\|b_t\|_{L^2}^2+C\|b\|_{L^3}^\frac{2}{3}\|\nabla^2 b\|_{L^2}^\frac{4}{3}\|\nabla u\|_{L^2}^2\nonumber\\
&\le \frac{1}{2}\|b_t\|_{L^2}^2+C\|b\|_{L^3}^\frac{1}{2}\|\nabla^2 b\|_{L^2}^2+C\|b\|_{L^3}\|\nabla u\|_{L^2}^6\nonumber\\
&\le \frac{1}{2}\|b_t\|_{L^2}^2+CM_0^\frac{1}{2}\|\nabla^2 b\|_{L^2}^2+CM_0\|\nabla u\|_{L^2}^2.
\end{align*}
Noticing that $\eqref{a1}_3$ and $p=\rho\theta$ imply that
\begin{align}\label{ky}
p_t=-{\rm div}(pu)-\rho\theta{\rm div}u+\mu(\nabla u+(\nabla u)^\top):\nabla u+\lambda(\divv u)^2+\Delta\theta
+|\curl b|^2.
\end{align}
Substituting \eqref{ky} into $J_4$, and using H\"older's, Young's, and Gagliardo-Nirenberg inequalities, \eqref{xss},
 and integration
by parts, one obtains
\begin{align}\label{yq}
J_4&=-\frac{1}{2\mu+\lambda}\int pu\cdot\nabla Fdx+\frac{1}{2\mu+\lambda}\int \rho\theta{\rm div}uFdx\nonumber\\
&\quad+\frac{\mu}{2\mu+\lambda}\int(\nabla u+(\nabla u)^\top):(\nabla F\otimes u)dx+\frac{\lambda}{2\mu+\lambda}\int{\rm div}uu\cdot\nabla Fdx\nonumber\\
&\quad+\frac{1}{2\mu+\lambda}\int(\mu\Delta u+(\mu+\lambda)\nabla{\rm div}u)\cdot uFdx+\frac{1}{2\mu+\lambda}\int\nabla\theta\cdot\nabla Fdx
\nonumber\\
&\quad+\frac{1}{2\mu+\lambda}\int|\curl b|^2Fdx\nonumber\\
&=-\frac{2}{2\mu+\lambda}\int pu\cdot\nabla Fdx+\frac{\mu}{2\mu+\lambda}\int(\nabla u+(\nabla u)^\top):(\nabla F\otimes u)dx\nonumber\\
&\quad+\frac{\lambda}{2\mu+\lambda}\int{\rm div}uu\cdot\nabla Fdx+\frac{1}{2\mu+\lambda}\int\nabla\theta\cdot\nabla Fdx
+\frac{1}{2\mu+\lambda}\int\rho u_t\cdot uFdx\nonumber\\
&\quad+\frac{1}{2\mu+\lambda}\int\rho u\cdot\nabla u\cdot uFdx+\frac{1}{2\mu+\lambda}\int b\otimes b:\nabla(uF)dx\nonumber\\
&\quad-\frac{1}{2(2\mu+\lambda)}
\int|b|^2{\rm div}(uF)dx+\frac{1}{2\mu+\lambda}
\int|\curl b|^2Fdx\nonumber\\
&\le C\|\nabla F\|_{L^2}(\|\rho u\theta\|_{L^2}+\||u||\nabla u|\|_{L^2}+\|\nabla\theta\|_{L^2}
+\||u||b|^2\|_{L^2})+\frac{1}{12}\int\rho|u_t|^2dx\nonumber\\
&\quad+C\int\rho|u|^2|F|^2dx+C\bar{\rho}\||u||\nabla u|\|_{L^2}^2
+C\|\nabla u\|_{L^2}\|b\|_{L^\infty}\|b\|_{L^3}\|F\|_{L^6}\nonumber\\
&\quad+C\|\nabla F\|_{L^2}\|b\|_{L^3}\|\nabla b\|_{L^6}+C\|F\|_{L^6}\|\nabla^2 b\|_{L^2}\|b\|_{L^3}\nonumber\\
&\le (CM_0^\frac{1}{4}+C\bar{\rho}+C\bar{\rho}^\frac{4}{3})\||u||\nabla u|\|_{L^2}^2+\frac{1}{12}\|\sqrt{\rho}u_t\|_{L^2}^2
+CM_0^{-\frac{1}{4}}\|\nabla F\|_{L^2}^2+CM_0^\frac{1}{4}\|\nabla\theta\|_{L^2}^2
\nonumber\\
&\quad+C\|\nabla u\|_{L^2}^2\|b\|_{L^3}^\frac{8}{3}\|\nabla^2 b\|_{L^2}^\frac{4}{3}
+C\|b\|_{L^3}^2\|\nabla^2 b\|_{L^2}^2
+C\bar{\rho}^\frac{4}{3}\|\rho\|_{L^1}^\frac{1}{3}\|\sqrt{\rho}\theta\|_{L^2}^2\|\nabla\theta\|_{L^2}^2
\nonumber\\
&\le CM_0^\frac{1}{4}\||u||\nabla u|\|_{L^2}^2+CM_0^{-\frac{1}{4}}\|\nabla F\|_{L^2}^2
+CM_0^\frac{1}{4}\|\nabla\theta\|_{L^2}^2+CM_0^4\|\nabla u\|_{L^2}^6\nonumber\\
&\quad+CM_0^2\|\nabla^2 b\|_{L^2}^2+CM_0^\frac{8}{3}\|\nabla\theta\|_{L^2}^2
+\frac{1}{12}\|\sqrt{\rho}u_t\|_{L^2}^2.
\end{align}
Taking the operator $\divv$ on both side of $\eqref{a1}_2$ gives rise to
\begin{align}\label{pss}
\Delta F={\rm div}(\rho u_t+\rho u\cdot\nabla u+b\cdot\nabla b),
\end{align}
which together with the standard elliptic estimates yields
\begin{align}\label{yyu}
\|\nabla F\|_{L^2}&\le C\bar{\rho}^\frac{1}{2}\|\sqrt{\rho}u_t\|_{L^2}
+C\bar{\rho}\||u||\nabla u|\|_{L^2}+C\|b\|_{L^3}\|\nabla b\|_{L^6}\nonumber\\
&\le CM_0\|\sqrt{\rho}u_t\|_{L^2}+CM_0\||u||\nabla u|\|_{L^2}+CM_0\|\nabla^2 b\|_{L^2}.
\end{align}
Substituting \eqref{yyu} into \eqref{yq}, and using \eqref{ee1}, we have
\begin{align*}
J_4\le CM_0^\frac{1}{4}\||u||\nabla u|\|_{L^2}^2+\Big(C_1M_0^\frac{7}{4}
+\frac{1}{12}\Big)\|\sqrt{\rho}u_t\|_{L^2}^2
+CM_0^\frac{7}{4}\|\nabla^2 b\|_{L^2}^2+CM_0^4\|\nabla u\|_{L^2}^2.
\end{align*}
Similarly, putting \eqref{ky} into $J_5$, one obtains
\begin{align}
J_5&=-\frac{1}{2\mu+\lambda}\int pu\cdot\nabla |b|^2dx+\frac{\mu}{2(2\mu+\lambda)}\int(\nabla u+(\nabla u)^\top):(\nabla |b|^2\otimes u)dx\nonumber\\
&\quad+\frac{\lambda}{2(2\mu+\lambda)}\int{\rm div}uu\cdot\nabla |b|^2dx+\frac{1}{2(2\mu+\lambda)}\int\nabla\theta\cdot\nabla |b|^2dx
\nonumber\\
&\quad+\frac{1}{2(2\mu+\lambda)}\int\rho u\cdot\nabla u\cdot u|b|^2dx+\frac{1}{2(2\mu+\lambda)}\int b\otimes b:\nabla(u|b|^2)dx\nonumber\\
&\quad-\frac{1}{4(2\mu+\lambda)}
\int|b|^2{\rm div}(u|b|^2)dx+\frac{1}{2(2\mu+\lambda)}\int|{\rm rot}b|^2|b|^2dx\nonumber\\
&\quad+\frac{1}{2(2\mu+\lambda)}\int\rho u_t\cdot u|b|^2dx\nonumber\\
&\le C\||b||\nabla b|\|_{L^2}(\|\rho u\theta\|_{L^2}+\||u||\nabla u|\|_{L^2}+\|\nabla\theta\|_{L^2}
+\||u||b|^2\|_{L^2})+\frac{1}{12}\int\rho|u_t|^2dx\nonumber\\
&\quad+C\int\rho|u|^2|b|^4dx+C\bar{\rho}\||u||\nabla u|\|_{L^2}^2+C\|\nabla u\|_{L^2}\||b|^4\|_{L^2}
+C\|b\|_{L^3}^2\|\nabla b\|_{L^6}^2\nonumber\\
&\le CM_0^\frac{1}{4}\||u||\nabla u|\|_{L^2}^2+CM_0^{-\frac{1}{4}}\||b||\nabla b|\|_{L^2}^2
+CM_0^\frac{1}{4}\|\nabla\theta\|_{L^2}^2+C\|u\|_{L^6}^2\|b\|_{L^\infty}^2\|b\|_{L^3}^2\nonumber\\[4pt]
&\quad+C\|\nabla u\|_{L^2}\|b\|_{L^8}^4
+C\|b\|_{L^3}^2\|\nabla^2 b\|_{L^2}^2+\frac{1}{12}\int\rho|u_t|^2dx\nonumber\\
&\le  CM_0^\frac{1}{4}\||u||\nabla u|\|_{L^2}^2+CM_0^{-\frac{1}{4}}\|b\|_{L^3}^2\|\nabla^2 b\|_{L^2}^2
+CM_0^\frac{1}{4}\|\nabla\theta\|_{L^2}^2+C\|b\|_{L^3}^2\|\nabla^2 b\|_{L^2}^2\nonumber\\
&\quad+C\|\nabla u\|_{L^2}^2\|b\|_{L^3}^\frac{8}{3}\|\nabla^2 b\|_{L^2}^\frac{4}{3}+
C\|\nabla u\|_{L^2}\|b\|_{L^3}^\frac{7}{3}\|\nabla^2 b\|_{L^2}^\frac{5}{3}+\frac{1}{12}\int\rho|u_t|^2dx
\nonumber\\
&\le  CM_0^\frac{1}{4}\||u||\nabla u|\|_{L^2}^2+CM_0^{-\frac{1}{4}}\|b\|_{L^3}^2\|\nabla^2 b\|_{L^2}^2
+CM_0^\frac{1}{4}\|\nabla\theta\|_{L^2}^2+C\|b\|_{L^3}^2\|\nabla^2 b\|_{L^2}^2\nonumber\\
&\quad+C\|b\|_{L^3}^4\|\nabla u\|_{L^2}^6+C\|b\|_{L^3}^2\|\nabla^2 b\|_{L^2}^2+C\|b\|_{L^3}^7\|\nabla u\|_{L^2}^6
\nonumber\\
&\quad+C\|b\|_{L^3}^\frac{7}{5}\|\nabla^2 b\|_{L^2}^2
+\frac{1}{12}\int\rho|u_t|^2dx\nonumber\\
&\le  \frac{1}{12}\int\rho|u_t|^2dx+CM_0^\frac{1}{4}\||u||\nabla u|\|_{L^2}^2+CM_0^\frac{7}{5}\|\nabla^2 b\|_{L^2}^2
+CM_0^\frac{1}{4}\|\nabla\theta\|_{L^2}^2+CM_0^4\|\nabla u\|_{L^2}^6,
\end{align}
where we have used the following fact
\begin{align*}
\|b\|_{L^8}^4\le \|b\|_{L^3}^\frac{2}{3}\|b\|_{L^{12}}^\frac{10}{3}
\le C\|b\|_{L^3}^\frac{2}{3}\||b||\nabla b|\|_{L^2}^\frac{10}{3}
\le C\|b\|_{L^3}^\frac{7}{3}\|\nabla b\|_{L^6}^\frac{5}{3}\le C\|b\|_{L^3}^\frac{7}{3}\|\nabla^2 b\|_{L^2}^\frac{5}{3}.
\end{align*}
Using Young's inequality and \eqref{ee1}, we have
\begin{align*}
J_6&\le \frac{1}{12}\int\rho|u_t|^2dx+C\int\rho|u|^2|\nabla u|^2dx
\le \frac{1}{12}\|\sqrt{\rho}u_t\|_{L^2}^2+CM_0^2\||u||\nabla u|\|_{L^2}^2.
\end{align*}
Substituting the above estimates on $J_i\ (i=3, 4, 5, 6)$ into \eqref{mmy} yields
\begin{align}\label{uub}
&\frac{1}{2}\frac{d}{dt}\big(\mu\|\nabla u\|_{L^2}^2+(\mu+\lambda)\|{\rm div}u\|_{L^2}^2\big)
+\frac{1}{2}\|\sqrt{\rho}u_t\|_{L^2}^2\nonumber\\
&\le \frac{d}{dt}\int \Big(\frac{1}{2}|b|^2{\rm div}u-b\cdot\nabla u\cdot b+p{\rm div}u\Big)dx-\frac{1}{2(2\mu+\lambda)}\frac{d}{dt}\int p^2dx\nonumber\\
&\quad+CM_0^\frac{1}{4}\||u||\nabla u|\|_{L^2}^2+CM_0^\frac{1}{2}\|\nabla^2 b\|_{L^2}^2
+CM_0\|\nabla u\|_{L^2}^2,
\end{align}
provided $M_0\le \epsilon_1=\min\Big\{\epsilon_1, \Big(\frac{1}{4C_1}\Big)^\frac{4}{7}\Big\}$.
Integrating \eqref{uub} over $[0, T]$, and using Cauchy-Schwarz inequality, we have
\begin{align}
&\mu\|\nabla u\|_{L^2}^2+(\mu+\lambda)\|{\rm div}u\|_{L^2}^2+\int_0^T\|\sqrt{\rho}u_t\|_{L^2}^2dt\nonumber\\
&\le \mu\|\nabla u_0\|_{L^2}^2+(\mu+\lambda)\|{\rm div}u_0\|_{L^2}^2-2\int\rho_0\theta_0{\rm div}u_0dx
+\frac{1}{2\mu+\lambda}\int\rho_0^2\theta_0^2dx\nonumber\\
&\quad+C\|b_0\|_{L^3}\|b_0\|_{L^6}\|\nabla u_0\|_{L^2}+C\|b\|_{L^3}\|b\|_{L^6}\|\nabla u\|_{L^2}
+\frac{1}{\mu+\lambda}\int\rho^2\theta^2dx\nonumber\\
&\quad+(\mu+\lambda)\int|{\rm div}u|^2dx+CM_0^\frac{1}{4}\int_0^T\||u||\nabla u|\|_{L^2}^2dt
+CM_0^\frac{1}{2}\int_0^T\|\nabla^2 b\|_{L^2}^2dt\nonumber\\
&\quad+CM_0\int_0^T\|\nabla u\|_{L^2}^2dt\nonumber\\
&\le \mu\|\nabla u_0\|_{L^2}^2+(\mu+\lambda)\|{\rm div}u_0\|_{L^2}^2+\frac{2\mu+\lambda}{2}\|{\rm div}u_0\|_{L^2}^2
+\frac{3\bar{\rho}}{2\mu+\lambda}\|\sqrt{\rho_0}\theta_0\|_{L^2}^2\nonumber\\
&\quad+\frac{\mu}{2}\|\nabla u_0\|_{L^2}^2+\frac{C}{\mu}M_0^2\|\nabla b_0\|_{L^2}^2
+\frac{\mu}{2}\|\nabla u\|_{L^2}^2+\frac{C}{\mu}M_0^2\|\nabla b\|_{L^2}^2\nonumber\\
&\quad+\frac{CM_0}{\mu+\lambda}\|\sqrt{\rho}\theta\|_{L^2}^2+(\mu+\lambda)\|{\rm div}u\|_{L^2}^2
+CM_0^\frac{1}{4}\int_0^T\||u||\nabla u|\|_{L^2}^2dt\nonumber\\
&\quad+CM_0^\frac{1}{2}\int_0^T\|\nabla^2 b\|_{L^2}^2dt
+CM_0\int_0^T\|\nabla u\|_{L^2}^2dt,
\end{align}
which yields that
\begin{align}\label{dd1}
&\frac{1}{2}\|\nabla u\|_{L^2}^2+\frac{1}{\mu}\int_0^t\|\sqrt{\rho}u_t\|_{L^2}^2dt\nonumber\\
&\le \frac{3}{2}\|\nabla u_0\|_{L^2}^2+\frac{3(\mu+\lambda)}{\mu}\|\nabla u_0\|_{L^2}^2+\frac{3(2\mu+\lambda)}{2\mu}\|\nabla u_0\|_{L^2}^2
+\frac{3M_0^2}{\mu(2\mu+\lambda)}\|\sqrt{\rho_0}\theta_0\|_{L^2}^2\nonumber\\
&\quad+CM_0^2\|\nabla b_0\|_{L^2}^2+C_2M_0^2\|\nabla b\|_{L^2}^2
+CM_0\|\sqrt{\rho}\theta\|_{L^2}^2
+CM_0^\frac{1}{4}\int_0^T\||u||\nabla u|\|_{L^2}^2dt\nonumber\\
&\quad+C_2M_0^\frac{1}{2}\int_0^T\|\nabla^2 b\|_{L^2}^2dt
+CM_0\int_0^T\|\nabla u\|_{L^2}^2dt,
\end{align}
where we have used
\begin{align*}
\|\divv u_0\|_{L^2}^2\le 3\|\nabla u_0\|_{L^2}^2.
\end{align*}

It follows from $\eqref{a1}_4$ that
\begin{align}\label{op}
&\frac{d}{dt}\|\nabla b\|_{L^2}^2+\|b_t\|_{L^2}^2+\|\nabla^2 b\|_{L^2}^2\nonumber\\
&=\int|b_t-\Delta b|^2dx=\int|b\cdot\nabla u-u\cdot\nabla b-b{\rm div}u|^2dx\nonumber\\
&\le C\|\nabla u\|_{L^2}^2\|b\|_{L^\infty}^2+C\|u\|_{6}^2\|\nabla b\|_{L^3}^2
\le C\|\nabla u\|_{L^2}^2\|b\|_{L^3}^\frac{2}{3}\|\nabla^2 b\|_{L^2}^\frac{4}{3}\nonumber\\
&\le C\|b\|_{L^3}\|\nabla u\|_{L^2}^6+C\|b\|_{L^3}^\frac{1}{2}\|\nabla^2 b\|_{L^2}^2\nonumber\\
&\le CM_0\|\nabla u\|_{L^2}^2+CM_0^\frac{1}{2}\|\nabla^2 b\|_{L^2}^2.
\end{align}
Integrating \eqref{op} over $[0, T]$ leads to
\begin{align}\label{dd2}
&\sup_{0\le t\le T}\|\nabla b\|_{L^2}^2+\int_0^T(\|b_t\|_{L^2}^2+\|\nabla^2 b\|_{L^2}^2)dt\nonumber\\
&\le \|\nabla b_0\|_{L^2}^2+CM_0\int_0^T\|\nabla u\|_{L^2}^2dt
+C_3M_0^\frac{1}{2}\int_0^T\|\nabla^2 b\|_{L^2}^2dt.
\end{align}
Adding \eqref{dd2} to \eqref{dd1}, we get
\begin{align}\label{3.51}
&\frac{1}{2}\|\nabla u\|_{L^2}^2+\frac{1}{2}\|\nabla b\|_{L^2}^2+\int_0^T\Big(\frac{1}{\mu}\|\sqrt{\rho}u_t\|_{L^2}^2+\frac{1}{2}\|b_t\|_{L^2}^2
+\frac{1}{2}\|\nabla^2 b\|_{L^2}^2\Big)dt\nonumber\\
&\le \frac{15\mu+9\lambda}{2\mu}\|\nabla u_0\|_{L^2}^2+\|\nabla b_0\|_{L^2}^2
+\frac{3M_0^2}{\mu(2\mu+\lambda)}\|\sqrt{\rho_0}\theta_0\|_{L^2}^2
+CM_0^2\|\nabla b_0\|_{L^2}^2\nonumber\\
&\quad+CM_0\|\sqrt{\rho}\theta\|_{L^2}^2
+CM_0^\frac{1}{4}\int_0^T\||u||\nabla u|\|_{L^2}^2dt
+CM_0\int_0^T\|\nabla u\|_{L^2}^2dt,
\end{align}
provided $M_0\le \epsilon_2=\min\Big\{\epsilon_1, \Big(\frac{1}{4C_2}\Big)^2,
\Big(\frac{1}{4C_3}\Big)^2\Big\}$.
Hence, the desired \eqref{vxx} follows from \eqref{3.51}.
\hfill $\Box$

\begin{lemma}\label{l36}
Under the conditions of Proposition \ref{p1}, it holds that
\begin{align}\label{xxc}
&\sup_{0\le t\le T}\Big(\|\sqrt{\rho}u\|_{L^2}^2+\|\sqrt{\rho}E\|_{L^2}^2+\|b\|_{L^2}^2+
\|\nabla u\|_{L^2}^2+\|\nabla b\|_{L^2}^2\Big)\nonumber\\
&\quad+\int_0^T\Big(\frac{\mu}{2}\|\nabla u\|_{L^2}^2+\|\nabla b\|_{L^2}^2+\frac{1}{2}\|\nabla\theta\|_{L^2}^2
+\frac{c_1c_2}{2}\||u||\nabla u|\|_{L^2}^2\Big)dt\nonumber\\
&\quad+\int_0^T\Big(\frac{2}{\mu}\|\sqrt{\rho}u_t\|_{L^2}^2
+\|b_t\|_{L^2}^2+\frac{1}{2}\|\nabla^2 b\|_{L^2}^2\Big)dt
\le \frac{7}{4}\hbar K,
\end{align}
provided
\begin{align*}
M_0\le \epsilon_3=\min&\Big\{\epsilon_2, \sqrt{\frac{\mu(2\mu+\lambda)}{40}}, \sqrt{\frac{3}{20C_4}},
\Big(\frac{3}{20C_4}\Big)^\frac{3}{5}, \frac{3}{20\hbar C_4},
\Big(\frac{3}{20\hbar KC_4}\Big)^\frac{3}{5}, \sqrt{\frac{1}{2C_4}}, \nonumber\\
&\quad\frac{\mu}{2C_4},\Big(\frac{3c_1c_2-5}{6C_4}\Big)^4\Big\}.
\end{align*}
Here $c_2$ is an absolute constant and $c_1$ is the same as that of in Lemma
\ref{lem34}.
\end{lemma}
{\it Proof.}
Based on Lemmas \ref{l31}--\ref{l35}, and adding $\eqref{bb}+\eqref{vvq}+c_2\times\eqref{sk}+\eqref{vxx}$ altogether for enough
large constant $c_2$,
it follows from \eqref{ee1} that
\begin{align}\label{ss5}
&\sup_{0\le t\le T}\Big(\|\sqrt{\rho}u\|_{L^2}^2+\|\sqrt{\rho}E\|_{L^2}^2+\|b\|_{L^2}^2+
\|\nabla u\|_{L^2}^2+\|\nabla b\|_{L^2}^2\Big)\nonumber\\
&\quad+\int_0^T\Big(\mu\|\nabla u\|_{L^2}^2+\|\nabla b\|_{L^2}^2+\|\nabla\theta\|_{L^2}^2+c_1c_2\||u||\nabla u|\|_{L^2}^2\Big)dt\nonumber\\
&\quad+\int_0^T\Big(\frac{2}{\mu}\|\sqrt{\rho}u_t\|_{L^2}^2
+\|b_t\|_{L^2}^2+\|\nabla^2 b\|_{L^2}^2\Big)dt\nonumber\\
&\le \|\sqrt{\rho_0}E_0\|_{L^2}^2+\|\sqrt{\rho_0}u_0\|_{L^2}^2+\|b_0\|_{L^2}^2
+\frac{15\mu+9\lambda}{\mu}\|\nabla u_0\|_{L^2}^2+2\|\nabla b_0\|_{L^2}^2\nonumber\\
&\quad+\frac{6M_0^2}{\mu(2\mu+\lambda)}\|\sqrt{\rho_0}\theta_0\|_{L^2}^2
+CM_0^2\|\nabla b_0\|_{L^2}^2+CM_0^\frac{5}{3}\|\nabla u_0\|_{L^2}^2+CM_0\sup_{0\le t\le T}\|\sqrt{\rho}\theta\|_{L^2}^2\nonumber\\
&\quad+CM_0^\frac{8}{3}\int_0^T\|\nabla\theta\|_{L^2}^2dt
+C\sup_{0\le t\le T}\|\rho^\frac{1}{4}u\|_{L^4}^4+\Big(\frac{5}{6}+CM_0^\frac{1}{4}\Big)\int_0^T\||u||\nabla u|\|_{L^2}^2dt\nonumber\\
&\quad+CM_0^2\int_0^T\|\nabla^2 b\|_{L^2}^2dt+CM_0^3\int_0^T\|\nabla u\|_{L^2}^2dt\nonumber\\
&\le \hbar K+\frac{6M_0^2}{\mu(2\mu+\lambda)}K+C_4M_0^2K+C_4M_0^\frac{5}{3}K+C_4M_0\hbar K+C_4M_0^\frac{5}{3}\hbar^2K^2\nonumber\\
&\quad+C_4M_0^\frac{8}{3}\int_0^T\|\nabla\theta\|_{L^2}^2dt
+\Big(\frac{5}{6}+C_4M_0^\frac{1}{4}\Big)\int_0^T\||u||\nabla u|\|_{L^2}^2dt
+C_4M_0^2\int_0^T\|\nabla^2 b\|_{L^2}^2dt\nonumber\\
&\quad+C_4M_0\int_0^T\|\nabla u\|_{L^2}^2dt,
\end{align}
where we have used
\begin{align*}
\|\rho^\frac{1}{4}u\|_{L^4}^4&\le \|\rho\|_{L^\infty}^\frac{1}{2}\|\sqrt{\rho}u\|_{L^2}\|u\|_{L^6}^3
\le C\|\rho\|_{L^\infty}^\frac{5}{6}\|\rho\|_{L^1}^\frac{2}{3}\|\nabla u\|_{L^2}^4\le CM_0^\frac{5}{3}\hbar^2K^2.
\end{align*}
Thus, it follows from \eqref{ss5} that
\begin{align*}
&\sup_{0\le t\le T}(\|\sqrt{\rho}u\|_{L^2}^2+\|\sqrt{\rho}E\|_{L^2}^2+\|b\|_{L^2}^2+
\|\nabla u\|_{L^2}^2+\|\nabla b\|_{L^2}^2)\nonumber\\
&\quad+\int_0^T(\frac{\mu}{2}\|\nabla u\|_{L^2}^2+\|\nabla b\|_{L^2}^2+\frac{1}{2}\|\nabla\theta\|_{L^2}^2
+\frac{c_1c_2}{2}\||u||\nabla u|\|_{L^2}^2)dt\nonumber\\
&\quad+\int_0^T(\frac{2}{\mu}\|\sqrt{\rho}u_t\|_{L^2}^2
+\|b_t\|_{L^2}^2+\frac{1}{2}\|\nabla^2 b\|_{L^2}^2)dt\nonumber\\
&\le \hbar K+\frac{3}{4}\hbar K=\frac{7}{4}\hbar K,
\end{align*}
provided
\begin{align*}
M_0\le \epsilon_3=\min&\Big\{\epsilon_2, \sqrt{\frac{\mu(2\mu+\lambda)}{40}}, \sqrt{\frac{3}{20C_4}},
\Big(\frac{3}{20C_4}\Big)^\frac{3}{5}, \frac{3}{20\hbar C_4},
\Big(\frac{3}{20\hbar KC_4}\Big)^\frac{3}{5}, \sqrt{\frac{1}{2C_4}}, \nonumber\\
&\quad\frac{\mu}{2C_4},\Big(\frac{3c_1c_2-5}{6C_4}\Big)^4\Big\}.
\end{align*}
The proof of Lemma \ref{l36} is finished.
\hfill $\Box$

\begin{lemma}
Under the conditions of Proposition \ref{p1}, it holds that
\begin{align}\label{ttr}
\sup_{0\le t\le T}t\|\nabla b\|_{L^2}^2+\int_0^Tt\big(\|b_t\|_{L^2}^2+\|\nabla ^2 b\|_{L^2}^2\big)dt\le C.
\end{align}
\end{lemma}
{\it Proof.}
Using H\"older's and Gagliardo-Nirenberg inequalities, we have
\begin{align*}
\frac{d}{dt}\|\nabla b\|_{L^2}^2+\|b_t\|_{L^2}^2+\|\nabla^2 b\|_{L^2}^2
&=\int|b\cdot\nabla u-u\cdot\nabla b-b{\rm div}u|^2dx\nonumber\\
&\le C\|u\|_{L^6}^2\|\nabla b\|_{L^3}^2+C\|\nabla u\|_{L^2}^2\|b\|_{L^\infty}^2\nonumber\\
&\le C\|\nabla u\|_{L^2}^2\|\nabla b\|_{L^2}\|\nabla^2 b\|_{L^2}\nonumber\\
&\le \frac{1}{2}\|\nabla^2 b\|_{L^2}^2+C\|\nabla u\|_{L^2}^4\|\nabla b\|_{L^2}^2.
\end{align*}
which implies that
\begin{align*}
\frac{d}{dt}\Big(t\|\nabla b\|_{L^2}^2\Big)+t\|b_t\|_{L^2}^2+\frac{t}{2}\|\nabla^2 b\|_{L^2}^2
&\le \|\nabla b\|_{L^2}^2+C\|\nabla u\|_{L^2}^4(t\|\nabla b\|_{L^2}^2).
\end{align*}
This together with Gronwall's inequality and \eqref{xxc} leads to the desired \eqref{ttr}.
\hfill $\Box$

\begin{lemma}\label{lem38}
Under the conditions of Proposition \ref{p1}, it holds that
\begin{align}\label{ghh}
0\le \rho\le \frac{3\bar{\rho}}{2},
\end{align}
provided  $M_0\le \epsilon_4=\min\Big\{\epsilon_3, \frac{\big(\log\frac{3}{2}\big)^3}{C_5^3}\Big\}$.
\end{lemma}
{\it Proof.}
The first inequality of \eqref{ghh} is obvious (see \cite[p. 43]{F04}). We only need to prove the second inequality of \eqref{ghh}. To this end, motivated by \cite{D1997,L1998} (see also \cite{WZ17}),
for any given $(x, t)\in \Bbb R^3\times[0, T]$, denote
\begin{align}
\rho^\delta(y,s)=\rho(y, s)+\delta\exp\Big\{-\int_0^s\divv(X(\tau; x, t), \tau)d\tau\Big\}>0
\end{align}
where
$X(s; x, t)$ is given by
\begin{align}
\left\{
\begin{array}{ll}
\displaystyle
\frac{d}{ds}X(s; x, t)=u(X(s; x, t), s), \quad 0\le s<t,\\
X(t; x, t)=x.
\end{array}
\right.
\end{align}
Using the fact that $\frac{d}{ds}(f(X(s; x, t), s)=(f_s+u\cdot\nabla f)(X(s; x, t), s)$,
 it follows from $\eqref{a1}_1$ that
 \begin{align}\label{tb1}
 \frac{d}{ds}\big(\log(\rho^\delta(X(s; x, t), s)\big)=-\divv u(X(s; x, t), s),
 \end{align}
which leads to
\begin{align}
Y'(s)=g(s)+b'(s),
\end{align}
where
\begin{align}
&Y(s)=\log\rho^\delta(X(s; x, t), s), \quad g(s)=-\frac{p(X(s; x, t), s)}{2\mu+\lambda},\notag \\
&b(s)=-\frac{1}{2\mu+\lambda}\int_0^s\Big(\frac{1}{2}|b(X(\tau; x, t), \tau)|^2+F(X(\tau; x, t), \tau)\Big)d\tau, \label{3.58}
\end{align}
and $F=(2\mu+\lambda)\divv u-p-\frac{1}{2}|b|^2
=(2\mu+\lambda)\divv u-\rho\theta-\frac12|b|^2$.

Rewrite $\eqref{a1}_2$ as
\begin{align}\label{dxy}
\partial_t\big[\Delta^{-1}{\rm div}(\rho u)\big]
-(2\mu+\lambda)\divv u+p+\frac{1}{2}|b|^2=-\Delta^{-1}\divv \divv (\rho u\otimes u)
+\Delta^{-1}\divv \divv (b\otimes b),
\end{align}
which implies that
\begin{align}\label{3.60}
F(X(\tau; x, t), \tau)
&=-\big[(-\Delta)^{-1}\divv (\rho u)\big]_{\tau}-(-\Delta)^{-1}{\rm div}{\rm div}(\rho u\otimes u)
+(-\Delta)^{-1}{\rm div}{\rm div}(b\otimes b)\nonumber\\
&=-\big[(-\Delta)^{-1}{\rm div}(\rho u)\big]_{\tau}-u\cdot\nabla(-\Delta)^{-1}{\rm div}(\rho u)
+u\cdot\nabla(-\Delta)^{-1}{\rm div}(\rho u)\nonumber\\
&\quad-(-\Delta)^{-1}{\rm div}{\rm div}(\rho u\otimes u)
+(-\Delta)^{-1}{\rm div}{\rm div}(b\otimes b)\nonumber\\
&=-\frac{d}{d\tau}\big[(-\Delta)^{-1}{\rm div}(\rho u)\big]+[u_i, R_{ij}](\rho u_j)+(-\Delta)^{-1}{\rm div}{\rm div}(b\otimes b),
\end{align}
where $[u_i, R_{ij}]=u_iR_{ij}-R_{ij}u_i$, and $R_{ij}=\partial_i(-\Delta)^{-1}\partial_j$ is the Riesz transform on $\Bbb R^3$. Hence we derive from \eqref{3.58} and \eqref{3.60} that
\begin{align}\label{xyx}
b(t)-b(0)
&\le\frac{1}{2\mu+\lambda}\int_0^t\Big[\frac{d}{d\tau}\big[(-\Delta)^{-1}{\rm div}(\rho u)\big]-[u_i, R_{ij}](\rho u_j)
-(-\Delta)^{-1}{\rm div}{\rm div}(b\otimes b)\Big]d\tau\nonumber\\
&\quad+\frac{1}{2(2\mu+\lambda)}\int_0^t\|b\|_{L^\infty}^2d\tau\nonumber\\
&\le\frac{1}{2\mu+\lambda}(-\Delta)^{-1}{\rm div}(\rho u)-\frac{1}{2\mu+\lambda}(-\Delta)^{-1}{\rm div}(\rho_0u_0)
+\frac{1}{2\mu+\lambda}\int_0^t\|[u_i, R_{ij}](\rho u_j)\|_{L^\infty}d\tau\nonumber\\
&\quad+\frac{1}{2\mu+\lambda}\int_0^t\|(-\Delta)^{-1}{\rm div}{\rm div}(b\otimes b)\|_{L^\infty}d\tau
+\frac{1}{2(2\mu+\lambda)}\int_0^t\|b\|_{L^\infty}^2d\tau\nonumber\\
&\le \frac{1}{2\mu+\lambda}\|(-\Delta)^{-1}{\rm div}(\rho u)\|_{L^\infty}
+\frac{1}{2\mu+\lambda}\|(-\Delta)^{-1}{\rm div}(\rho_0u_0)\|_{L^\infty}\nonumber\\
&\quad+\frac{1}{2\mu+\lambda}\int_0^t\|[u_i, R_{ij}](\rho u_j)\|_{L^\infty}d\tau+
\frac{1}{2\mu+\lambda}\int_0^t\|(-\Delta)^{-1}{\rm div}{\rm div}(b\otimes b)\|_{L^\infty}d\tau\nonumber\\
&\quad+\frac{1}{2(2\mu+\lambda)}\int_0^t\|b\|_{L^\infty}^2d\tau=\sum_{i=1}^5Z_i.
\end{align}

By Gagliardo-Nirenberg, Sobolev's, Calder{\'o}n-Zygmund,
 and H\"older's inequalities, \eqref{abb}, and \eqref{ee1}, one obtains
\begin{align}\label{3.62}
Z_1&\le \frac{C}{2\mu+\lambda}\|(-\Delta)^{-1}{\rm div}(\rho u)\|_{L^6}^\frac{1}{3}
\|\nabla(-\Delta)^{-1}{\rm div}(\rho u)\|_{L^4}^\frac{2}{3}\nonumber\\
&\le C\|\rho u\|_{L^2}^\frac{1}{3}\|\rho u\|_{L^4}^\frac{2}{3}
\le C\|\rho\|_{L^3}^\frac{1}{3}\|u\|_{L^6}^\frac{1}{3}\|\rho\|_{L^{12}}^\frac{2}{3}\|u\|_{L^6}^\frac{2}{3}\nonumber\\
&\le C\|\rho\|_{L^\infty}^\frac{15}{18}\|\rho\|_{L^1}^\frac{1}{6}\|\nabla u\|_{L^2}
\le CM_0^\frac{15}{9}.
\end{align}
Similarly to \eqref{3.62}, we have
\begin{align}
Z_2&\le CM_0^\frac{15}{9}.
\end{align}
For $Z_3$, we deduce from Gagliardo-Nirenberg inequality and Calder{\'o}n-Zygmund inequality that
\begin{align}\label{uuy}
Z_3&\le \frac{C}{2\mu+\lambda}\int_0^t\|[u_i, R_{ij}](\rho u_j)\|_{L^3}^\frac15
\|\nabla[u_i, R_{ij}](\rho u_j)\|_{L^4}^\frac45d\tau\nonumber\\
&\le C\int_0^t\|u\|_{L^6}^\frac{1}{5}\|\rho u\|_{L^6}^\frac15\|\nabla u\|_{L^6}^\frac45\|\rho u\|_{L^{12}}^\frac45d\tau\nonumber\\
&\le C\int_0^t\|\rho\|_{L^\infty}\|u\|_{L^6}^\frac{1}{5}\|\nabla u\|_{L^6}^\frac45
\Big(\|u\|_{L^6}^\frac34\|\nabla u\|_{L^6}^\frac14\Big)^\frac45d\tau\nonumber\\
&\le C\int_0^t\bar{\rho}\|\nabla u\|_{L^2}\|\nabla u\|_{L^6}d\tau.
\end{align}
Denote $w=\curl u$, then we have (see e.g., \cite[Theorem 11.25]{FN2017})
\begin{align}\label{syk}
\|\nabla u\|_{L^6}\le C\|w\|_{L^6}+C\|\divv u\|_{L^6}.
\end{align}
Taking the operators $\divv$ and $\curl$ on both sides of $\eqref{a1}_2$ respectively, we get
\begin{align}
\left\{
\begin{array}{ll}
\displaystyle
\Delta F={\rm div}(\rho u_t+\rho u\cdot\nabla u)+{\rm div}{\rm div}(b\otimes b),\\
\mu\Delta w=\nabla\times(\rho u_t+\rho u\cdot\nabla u+{\rm div}(b\otimes b),
\end{array}
\right.\label{opp}
\end{align}
which together with the standard elliptic estimates implies that
\begin{align}\label{3.68}
\|\nabla w\|_{L^2}+\|\nabla F\|_{L^2}&\le C\bar{\rho}^\frac{1}{2}\|\sqrt{\rho} u_t\|_{L^2}+C\bar{\rho}\||u||\nabla u|\|_{L^2}
+C\|b\nabla b|\|_{L^2}\nonumber\\
&\le C\bar{\rho}^\frac{1}{2}\|\sqrt{\rho} u_t\|_{L^2}+C\bar{\rho}\||u||\nabla u|\|_{L^2}+C\|b\|_{L^3}\|\nabla b\|_{L^6}\nonumber\\
&\le CM_0(\|\sqrt{\rho}u_t\|_{L^2}+\||u||\nabla u|\|_{L^2}+\|\nabla^2 b\|_{L^2}).
\end{align}
Substituting \eqref{syk} and \eqref{opp} into \eqref{uuy}, we infer from \eqref{3.68} and \eqref{xxc} that
\begin{align*}
Z_3&\le CM_0^2\int_0^t\|\nabla u\|_{L^2}(\|\nabla w\|_{L^6}+\|{\rm div}u\|_{L^6})d\tau\nonumber\\
&\le CM_0^2\int_0^t\|\nabla u\|_{L^2}\Big(\|\nabla w\|_{L^2}+\frac{1}{2\mu+\lambda}\|F\|_{L^6}
+\frac{1}{2\mu+\lambda}\|\rho\theta\|_{L^6}\Big)d\tau\nonumber\\
&\le CM_0^3\int_0^t\big(\|\nabla u\|_{L^2}^2+\|\nabla w\|_{L^2}^2+\|\nabla F\|_{L^2}^2+\|\nabla\theta\|_{L^2}^2\big)d\tau\nonumber\\
&\le CM_0^3\int_0^t\big(\|\nabla u\|_{L^2}^2+\|\sqrt{\rho}u_t\|_{L^2}^2+\||u||\nabla u|\|_{L^2}^2+\|\nabla^2 b\|_{L^2}^2+\|\nabla\theta\|_{L^2}^2\big)d\tau\nonumber\\
&\le CM_0^3.
\end{align*}
For $Z_4$, by H\"older's and Gagliardo-Nirenberg inequalities, \eqref{ee1}, \eqref{xxc}, and \eqref{ttr},  we have
\begin{align}\label{3.69}
Z_4&\le \frac{1}{2\mu+\lambda}\int_0^1\|(-\Delta)^{-1}{\rm div}{\rm div}(b\otimes b)\|_{L^\infty}d\tau
+\frac{1}{2\mu+\lambda}\int_1^t\|(-\Delta)^{-1}{\rm div}{\rm div}(b\otimes b)\|_{L^\infty}d\tau\nonumber\\
&\le \frac{C}{2\mu+\lambda}\int_0^1\|b\|_{L^3}^\frac{2}{3}\|\nabla^2 b\|_{L^2}^\frac{4}{3}d\tau
+\frac{C}{2\mu+\lambda}\int_1^t\|b\|_{L^3}^\frac{1}{2}\|\nabla^2 b\|_{L^2}^\frac{7}{6}\|\nabla b\|_{L^2}^\frac{1}{2}d\tau\nonumber\\
&\le CM_0^\frac{1}{3}\sup_{1\le \tau\le t}\|\nabla b\|_{L^2}^\frac{1}{2}\Big(\int_1^tt^{-\frac{7}{12}\cdot\frac{12}{5}}d\tau\Big)^\frac{5}{12}
\Big(\int_1^t\tau\|\nabla^2 b\|_{L^2}^2d\tau\Big)^\frac{7}{12}
+CM_0^\frac{2}{3}\Big(\int_0^1\|\nabla^2 b\|_{L^2}^2d\tau\Big)^\frac{2}{3}\nonumber\\
&\le CM_0^\frac{2}{3}+CM_0^\frac{1}{3}\le CM_0^\frac{1}{3}.
\end{align}
Here we have used the following Gagliardo-Nirenberg inequality
\begin{align*}
\|b\|_{L^\infty}\le C\|b\|_{L^3}^\frac{1}{3}\|\nabla^2 b\|_{L^2}^\frac{2}{3},
\quad\|b\|_{L^\infty}\le C\|\nabla b\|_{L^2}^\frac{1}{2}\|\nabla^2 b\|_{L^2}^\frac{1}{2}.
\end{align*}
Similarly to \eqref{3.69}, we have
\begin{align*}
Z_5&\le CM_0^\frac{1}{3}.
\end{align*}
Substituting the above estimates for $Z_i\ (i=1, 2, 3, 4, 5)$ into \eqref{xxc} yields
\begin{align}
b(t)-b(0)&\le CM_0^3+CM_0^\frac{1}{3}+CM_0^\frac{15}{9}\le C_5M_0^\frac{1}{3}\le \log\frac{3}{2},
\end{align}
provided $M_0\le \epsilon_4=\min\Big\{\epsilon_3, \frac{\big(\log\frac{3}{2}\big)^3}{C_5^3}\Big\}$.

Integrating \eqref{tb1} w.r.t.  $s$ over $[0, t]$, we get
\begin{align*}
\log\rho^\delta(x, t)&=\log[\rho_0(X(t; x, 0))+\delta]+\int_0^tg(\tau)d\tau+b(t)-b(0)\nonumber\\
&\le \log(\bar{\rho}+\delta)+\log\frac{3}{2}.
\end{align*}
Let $\delta\rightarrow 0^+$, we have
\begin{align*}
\rho\le \frac{3\bar{\rho}}{2}.
\end{align*}
This finishes the proof of Lemma \ref{lem38}. \hfill $\Box$

\begin{lemma}\label{l39}
Under the conditions of Proposition \ref{p1}, it holds that
\begin{align}
\sup_{0\le t\le T}\|b\|_{L^3}\le \frac{3}{2}M_0,
\end{align}
provided $M_0\le \epsilon_0=\min\Big\{\epsilon_4, \frac{3}{2C_6}\Big\}$.
\end{lemma}
{\it Proof.}
Multiplying $\eqref{a1}_4$ by $3|b|b$ and integrating by parts over $\Bbb R^3$, we have
\begin{align*}
\frac{d}{dt}\|b\|_{L^3}^3+3\int|b||\nabla b|^2dx+3\int|b||\nabla|b||^2dx
\le \int|b||\nabla b|^2dx+C\|\nabla u\|_{L^2}^2\|b\|_{L^\frac{9}{2}}^3.
\end{align*}
Consequently,
\begin{align}\label{ugg}
\frac{d}{dt}\|b\|_{L^3}^3+2\int|b||\nabla b|^2dx+3\int|b||\nabla|b||^2dx
\le C\|\nabla u\|_{L^2}^2\|b\|_{L^\frac{9}{2}}^3.
\end{align}
To deal with the right-hand side of \eqref{ugg}, we need to use the following variant of the Kato inequality
\begin{align*}
|\nabla|b|^\frac{3}{2}|=\frac{3}{2}|b|^\frac{1}{2}|\nabla|b||\le \frac{3}{2}|b|^\frac{1}{2}|\nabla b|,
\end{align*}
which combined with Sobolev's inequality and Galiardo-Nirenberg inequality leads to
\begin{align}\label{ssw}
\|b\|_{L^\frac{9}{2}}^3\le \|b\|_{L^3}^\frac32\|b\|_{L^9}^\frac32
= \|b\|_{L^3}^\frac32\||b|^\frac32\|_{L^6}
\le C\|b\|_{L^3}^\frac32\|\nabla(|b|^\frac32)\|_{L^2}
\le C\|b\|_{L^3}^\frac32\||b|^\frac12|\nabla b|\|_{L^2}.
\end{align}
Thus, putting \eqref{ssw} into \eqref{ugg}, we obtain from Cauchy-Schwarz inequality that
\begin{align*}
\frac{d}{dt}\|b\|_{L^3}^3+\int|b||\nabla b|^2dx\le C\|\nabla u\|_{L^2}^4\|b\|_{L^3}^3.
\end{align*}
This together with \eqref{xxc} and Gronwall's inequality yields
\begin{align*}
\sup_{0\le t\le T}\|b\|_{L^3}
\le \exp\Big\{C\int_0^T\|\nabla u\|_{L^2}^4dt\Big\}^\frac{1}{3}\|b_0\|_{L^3}\le C_6M_0^2\le \frac{3M_0}{2},
\end{align*}
provided $M_0\le \epsilon_0=\min\Big\{\epsilon_4, \frac{3}{2C_6}\Big\}$.
The lemma is completed.
\hfill $\Box$

Now, Proposition \ref{p1} is a direct consequence of Lemmas \ref{l31}--\ref{l39}.

{\it Proof of Proposition \ref{p1}.}
Define
\begin{align*}
T^{\#}:=\max\Big\{T'\in (0, T]\Big| \sup_{0\le t\le T'}\|\rho\|_{L^\infty}\le 2\bar{\rho},~~
\psi_{T'}\le 2\hbar K, ~\sup_{0\le t\le T'}\|b\|_{L^3}\le 2M_0\Big\}.
\end{align*}
Then, by Lemmas \ref{l31}--\ref{l39}, we have
\begin{align}\label{fhf}
\sup_{0\le t\le T}\|\rho\|_{L^\infty}\le \frac{3}{2}\bar{\rho}, \quad \psi_T\le \frac{7}{4}\hbar K,
\quad \sup_{0\le t\le T}\|b\|_{L^3}\le \frac{3}{2}M_0,\quad \forall T'\in(0, T^{\#}).
\end{align}
as long as $M_0\le \epsilon_0$ is sufficiently small.

If $T^{\#}<T$, noticing that $\psi_{T'}$, $\sup\limits_{0\le t\le T'}\|b\|_{L^3}$, and $\sup\limits_{0\le t\le T'}\|\rho\|_{L^\infty}$
are continuous on $[0, T]$, there is another time $T^{\#\#}\in (T^{\#}, T]$ such that
\begin{align*}
\sup_{0\le t\le T^{\#\#}}\|\rho\|_{L^\infty}\le \frac{3}{2}\bar{\rho}, \quad \psi_{T^{\#\#}}\le \frac{7}{4}\hbar K,
\quad \sup_{0\le t\le T^{\#\#}}\|b\|_{L^3}\le \frac{3}{2}M_0,
\end{align*}
which contradicts to the definition of $T^{\#}$. Thus, we have $T^{\#}=T$, and the conclusion
follows from \eqref{fhf} and the continuity of
$\psi_{T'}$, $\sup\limits_{0\le t\le T'}\|b\|_{L^3}$, and $\sup\limits_{0\le t\le T'}\|\rho\|_{L^\infty}$.
\hfill $\Box$

The following corollary is a straightforward consequence of Proposition \ref{p1} and Lemma \ref{l36}.
\begin{corollary}\label{qww}
Assume that $3\mu>\lambda$, and let the conditions in Proposition \ref{p1} be in force. Then there is a positive
constant $C$ depending only on $\mu$, $\lambda$, $\|\rho_0\|_{L^1}$, and $K$ such that
\begin{align*}
&\sup_{0\le t\le T}\Big(\|\rho\|_{L^\infty}+\|\sqrt{\rho}u\|_{L^2}^2+\|\sqrt{\rho}E\|_{L^2}^2+
\|\nabla u\|_{L^2}^2+\|b\|_{H^1}^2\Big)\nonumber\\
&\quad+\int_0^T\Big(\|\nabla u\|_{L^2}^2+\|\nabla b\|_{H^1}^2+\|\nabla\theta\|_{L^2}^2
+\||u||\nabla u|\|_{L^2}^2+\|\sqrt{\rho}u_t\|_{L^2}^2
+\|b_t\|_{L^2}^2\Big)dt
\le C,
\end{align*}
provided that $M_0\le \epsilon_0$.
\end{corollary}

\section{Proof of Theorem \ref{thm1}}

Let $\epsilon_0$ be the constant stated in Proposition \ref{p1} and suppose that the initial
data $(\rho_0, u_0, \theta_0, b_0)$ satisfies \eqref{qq} and \eqref{1.6}, and
\begin{align*}
 M_0\le \epsilon_0.
\end{align*}
According  to Lemma \ref{l22}, there is a unique local strong solution $(\rho, u, \theta, b)$ to the problem \eqref{a1}--\eqref{a2}. Let $T_{\max}$ be the maximal existence time to the solution. We will show that $T_{\max}=\infty$.
Suppose, by contradiction, that $T_{\max}<\infty$. Then, by virtue of Lemma \ref{pyy}, there holds
\begin{align}\label{4.1}
\lim_{T\rightarrow T_{\max}}\big(\|\rho\|_{L^\infty(0, T; L^\infty)}+\|u\|_{L^4(0, T; L^6)}\big)=\infty.
\end{align}
By Corollary \ref{qww}, for any $T\in (0, T_{\max})$, there exists a positive constant $\bar{C}$ independent of $T$ such that
\begin{align}\label{4.2}
\sup\limits_{0\le t\le T}\left(\|\rho\|_{L^\infty}+\|\nabla u\|_{L^2}^2\right)\le \bar{C},
\end{align}
which combined with Sobolev's inequality $\|u\|_{L^6}\leq C\|\nabla u\|_{L^2}$ gives
\begin{align}\label{4.3}
\int_0^{T_{\max}}\|u\|_{L^6}^4dt
\leq C\int_0^{T_{\max}}\|\nabla u\|_{L^2}^4dt\leq C\bar{C}^2T_{\max}<\infty.
\end{align}
From \eqref{4.2} and \eqref{4.3}, we derive that
\begin{align*}
\lim_{T\rightarrow T_{\max}}\big(\|\rho\|_{L^\infty(0, T; L^\infty)}+\|u\|_{L^4(0, T; L^6)}\big)<\infty,
\end{align*}
contradicting to \eqref{4.1}. This contradiction provides us that $T_{\max}=\infty$, and thus we obtain the
global strong solution. This finishes the proof of Theorem \ref{thm1}.
\hfill $\Box$

\end{document}